\documentclass[preprint,12pt]{elsarticle}

\usepackage{xcolor}

\usepackage{graphicx}
\usepackage{hyperref}

\usepackage{amssymb}
\usepackage{amsfonts}
\usepackage{amsmath}
\usepackage{amsthm}
\usepackage{verbatim}

\usepackage{bm}

\usepackage{multirow,multicol}
\usepackage{float}

\usepackage{tikz}
\usepackage{pgfplots}
\pgfplotsset{compat=1.18}
\usepackage{tikz-3dplot}

\usepackage{booktabs}
\usepackage{subcaption}
\usepackage{placeins}
\usepackage{doi}

\usepackage[ruled,vlined]{algorithm2e}

\usetikzlibrary{arrows.meta}
\usetikzlibrary{shapes.geometric,patterns,hobby}
\usetikzlibrary{spy}
\usetikzlibrary{shadings}
\usetikzlibrary{3d}
\usetikzlibrary{calc,decorations.pathreplacing}
\usetikzlibrary{decorations.pathmorphing}

\allowdisplaybreaks

\hypersetup{
	colorlinks=true,
	linkcolor=blue,
	citecolor=blue,
	urlcolor=blue
}

\begin{document}
	
	\begin{frontmatter}
		
		\title{Adaptive Time Windows for Discrete Adjoint Topology Optimization of Unsteady Flows}
		
		\author[aff1]{Zongyuan Liu}
		\ead{51275500058@stu.ecnu.edu.cn}
		
		\author[aff2]{Kentaro Yaji}
		\ead{yaji@mech.eng.osaka-u.ac.jp}
		
		\author[aff2]{Musaddiq Al Ali}
		\ead{drmusaddiqalali@gmail.com}
		
		\author[aff1]{Shengfeng Zhu\corref{cor1}}
		\ead{sfzhu@math.ecnu.edu.cn}
		
		\cortext[cor1]{Corresponding author}
		
		\affiliation[aff1]{
			organization={School of Mathematical Sciences,
				East China Normal University},
			addressline={500 Dongchuan Road},
			city={Shanghai},
			postcode={200241},
			country={China}
		}
		
		\affiliation[aff2]{
			organization={Department of Mechanical Engineering,
				Graduate School of Engineering,
				The University of Osaka},
			addressline={2-1 Yamadaoka},
			city={Suita},
			state={Osaka},
			postcode={565-0871},
			country={Japan}
		}
		
		\begin{abstract}
			Rather than prescribing an evaluation interval a priori,
			the proposed framework characterizes each evolving unsteady
			flow using a sequence of time windows. A consecutive-window
			convergence criterion is introduced to automatically identify
			a representative time window within its fully developed stage.
			The objective evaluation and discrete adjoint analysis are then
			carried out consistently over the identified representative
			time window. The framework is implemented using a regularized
			lattice Boltzmann method-based large-eddy simulation (LBM-LES)
			solver together with a partial bounce-back fluid-solid model.
			The consecutive-window convergence criterion is first validated
			using the backward-facing step flow. Cylinder-flow applications
			are then employed to investigate the influence of different flow
			regimes on the proposed framework. The wake-flow recovery problem
			verifies its effectiveness for unsteady topology optimization,
			while U-bend optimization further demonstrates its capability
			to identify and reorganize complex vortical structures.
		\end{abstract}
		
		
	\end{frontmatter}
	
	\section{Introduction}\label{intro}
	The design of fluid-flow systems is central to applications such as heat exchangers~\cite{Mohammed2011}, fluid-distribution networks~\cite{Chen2009Manifold}, aerodynamic components, chemical reactors, hydraulic machinery, and microfluidic devices~\cite{Zhang2023}. In these systems, geometry greatly influences pressure loss, flow uniformity, mixing, heat transfer, and drag. When the flow becomes unsteady, it also affects the mean flow and the evolution of coherent vortical structures~\cite{Jimenez2018}, making the characterization of representative flow behavior essential for design.
	
	Topology optimization provides a powerful framework for fluid-system design~\cite{Kobayashi2021,Hoghoj2020,Wang2018,Gaymann2019,Koga2013,Andreasen2009,Liu2026, Li2026a}. Unlike conventional shape optimization~\cite{Li2022,Li2026b}, it uses a spatially distributed design variable to distinguish fluid and solid regions. Darcy-type models~\cite{Borrvall2003,Wiker2007} typically represent intermediate material states, enabling the automatic creation, removal, merging, and separation of flow passages. Most existing fluid topology optimization studies, however, have focused on steady incompressible flows.
	
	For unsteady flow topology optimization, the lattice Boltzmann method (LBM) is particularly suitable because its explicit time-marching formulation naturally tracks flow evolution~\cite{Chen2017,Norgaard2016}. In this work, a regularized
	lattice Boltzmann large-eddy simulation (LBM-LES) solver incorporating the Smagorinsky subgrid-scale model~\cite{Smagorinsky1963} is adopted to resolve large-scale vortical structures while modeling unresolved scales.
	
	However, steady and unsteady flows exhibit fundamentally different long-time behaviors. From a dynamical-systems perspective, flow evolution traces a trajectory in phase space: a steady flow approaches a fixed point, whereas an unsteady flow may approach a limit cycle or, for more complex dynamics, a chaotic attractor~\cite{Strogatz2018,Guckenheimer1983}. Conventional steady-flow topology optimization evaluates the objective function and sensitivity at a converged stationary state. Because an unsteady trajectory does not converge to a unique state, the steady-state framework cannot be applied directly; instead, a representative portion of the trajectory must be identified for objective and sensitivity evaluation.
	
	Several approaches have been developed for this purpose. Nørgaard et al.~\cite{Norgaard2016} used a time-resolved discrete adjoint over a prescribed interval, requiring the primal trajectory for backward computation. Chen et al.~\cite{Chen2017} reduced the storage requirement by dividing the time horizon into subintervals and approximating the global sensitivity from local sensitivities. Nevertheless, the intervals remain prescribed and may become unrepresentative as topology updates alter the transient duration and vortex development. The influence of the selected averaging length has also been investigated~\cite{Putranto2026}, while massively parallel time integration has been used to address the computational cost~\cite{Katsumata2025}.
	To reduce the storage required by time-resolved adjoint computations, Ito et al.~\cite{Ito2026MeanFlow} constructed the adjoint equations from time-averaged primal states for sensitivity analysis of unsteady porous-cylinder flows simulated with the  homogenized lattice Boltzmann method (HLBM). In their framework, averaging begins after a prescribed spin-up period and continues until fluctuations of the averaged quantities fall below a specified tolerance. However, steady and unsteady flows are treated through separate computational branches, requiring prior identification of the flow regime. This requirement may become restrictive in topology optimization because structural updates can alter flow stability between design iterations.
	A separate approach constructs an auxiliary steady optimization problem from time-averaged flow quantities and estimated turbulence parameters~\cite{Ramalingom2025}. The resulting optimization therefore depends on the assumptions of the auxiliary model. These studies motivate a criterion that identifies a representative time window for each evolving design without prior specification of the resulting flow regime.
	
	The circular-cylinder wake provides a physical basis for such a criterion. It is a canonical example of the transition from a steady wake to periodic vortex shedding~\cite{Williamson1996}. After the initial transient, the instantaneous wake remains oscillatory, whereas its dominant frequency and time-averaged quantities approach stable values. Theoretical analyses likewise characterize the wake using the developed periodic state~\cite{Dusek1994}. Motivated by these observations, we identify the fully developed behavior through the convergence of averaged quantities over consecutive time windows. Unlike prescribed-window approaches, the representative window is determined independently and adaptively for each evolving structure, accommodating variations in flow-development duration and flow regime. Once the convergence criterion is satisfied, the same window is used for objective evaluation and backward discrete adjoint analysis, yielding a unified topology optimization framework for steady and unsteady flows.
	
	The paper is organized as follows. Section~\ref{Flow_solver_and_fluid-solid_modeling} presents the flow solver and fluid-solid model, Section~\ref{Statistical_adjoint_framework} introduces the proposed optimization framework and algorithm, Section~\ref{Numerical_Results} presents the validation and optimization results, and Section~\ref{Conclusion} concludes the paper.
	
	\section{Flow solver and fluid-solid modeling}\label{Flow_solver_and_fluid-solid_modeling}
	
	\subsection{Flow solver}
	To simulate two-dimensional (2D) unsteady incompressible flows, a regularized LBM-LES solver~\cite{Latt2006}
	is employed. The lattice Boltzmann formulation, LES treatment, boundary
	conditions, and fluid-solid interaction model are described in the following
	subsections.
	\subsubsection{Lattice Boltzmann method}
	\begin{figure}[htbp]
		\centering
		\begin{tikzpicture}[scale=1.0, >=stealth]
			
			\draw[gray!40] (-1,-1) grid (1,1);
			
			\node[circle, draw, fill=white, inner sep=2pt] (c0) at (0,0) {$c_0$};
			\node[circle, draw, fill=white, inner sep=2pt] (c1) at (1,0) {$c_1$};
			\node[circle, draw, fill=white, inner sep=2pt] (c2) at (0,1) {$c_2$};
			\node[circle, draw, fill=white, inner sep=2pt] (c3) at (-1,0) {$c_3$};
			\node[circle, draw, fill=white, inner sep=2pt] (c4) at (0,-1) {$c_4$};
			\node[circle, draw, fill=white, inner sep=2pt] (c5) at (1,1) {$c_5$};
			\node[circle, draw, fill=white, inner sep=2pt] (c6) at (-1,1) {$c_6$};
			\node[circle, draw, fill=white, inner sep=2pt] (c7) at (-1,-1) {$c_7$};
			\node[circle, draw, fill=white, inner sep=2pt] (c8) at (1,-1) {$c_8$};
			
			\draw[->, thick] (c0) -- (c1);
			\draw[->, thick] (c0) -- (c2);
			\draw[->, thick] (c0) -- (c3);
			\draw[->, thick] (c0) -- (c4);
			\draw[->, thick] (c0) -- (c5);
			\draw[->, thick] (c0) -- (c6);
			\draw[->, thick] (c0) -- (c7);
			\draw[->, thick] (c0) -- (c8);
			
		\end{tikzpicture}
		\caption{The D2Q9 lattice model.}
		\label{d2q9}
	\end{figure}
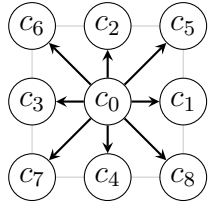
	Instead of directly solving the macroscopic Navier--Stokes equations, the LBM models fluid flows through the evolution of particle distribution functions on a discrete velocity lattice. The evolution of the distribution functions is governed by the Bhatnagar--Gross--Krook (BGK) lattice Boltzmann equation~\cite{Qian1992}:
	\begin{equation}\label{LBE}
		f_\alpha(\mathbf{x}_i+\mathbf{c}_\alpha\Delta t,t+\Delta t)=f_\alpha(\mathbf{x}_i,t)-\frac{1}{\tau}[f_\alpha(\mathbf{x}_i,t)-f_\alpha^{\rm eq}(\mathbf{x}_i,t)],
	\end{equation}
	where $\Delta t$ is the lattice time step, $\tau$ is the relaxation time, and $\mathbf{c}_\alpha$ denotes the discrete lattice velocity in the $\alpha$-th direction, with $\alpha=0,1,\ldots,8$. The function $f_\alpha(\mathbf{x},t)$ represents the particle distribution function associated with $\mathbf{c}_\alpha$ at position $\mathbf{x}_i$ and time $t$, while $f_\alpha^{\rm eq}$ is the corresponding equilibrium distribution function. In this study, the D2Q9 model~\cite{Qian1992} is adopted, which consists of one rest velocity, four axial velocities, and four diagonal velocities, as illustrated in Fig.~\ref{d2q9}.
	
	The discrete velocities of the D2Q9 model are given by
	\begin{equation}
		\begin{aligned}
			&\left[
			\mathbf c_0,\mathbf c_1,\mathbf c_2, \mathbf c_3, \mathbf c_4, \mathbf c_5, \mathbf c_6, \mathbf c_7,\mathbf c_8
			\right]\\
			&= \frac{\Delta x}{\Delta t}
			\begin{bmatrix}
				0 & 1 & 0 & -1 & 0 & 1 & -1 & -1 & 1\\
				0 & 0 & 1 & 0 & -1 & 1 & 1 & -1 & -1
			\end{bmatrix},
		\end{aligned}
	\end{equation}
	and for simplicity, lattice units are adopted throughout this study, i.e., $\Delta x=\Delta t=1$.
	
	For the D2Q9 lattice, the local equilibrium distribution function is expressed as
	\begin{equation}\label{feq}
		f_\alpha^{\rm eq}(\mathbf{x}_i,t)=w_\alpha\rho\bigg[1+\frac{\mathbf{c}_\alpha\cdot\mathbf{u}}{c_s^2}+\frac{(\mathbf{c}_\alpha\cdot\mathbf{u})^2}{2c_s^4}-\frac{|\mathbf{u}|^2}{2c_s^2}\bigg],
	\end{equation}
	where $w_\alpha$ is the weighting coefficient, $\rho$ is the fluid density, $\mathbf{u}$ is the macroscopic velocity, and $c_s$ denotes the lattice speed of sound. 
	
	The weighting coefficients are given by
	\begin{equation}
		w_\alpha=
		\begin{cases}
			4/9,  & \alpha=0,\\
			1/9,  & \alpha\in\{1,2,3,4\},\\
			1/36, & \alpha\in\{5,6,7,8\}.
		\end{cases}
	\end{equation}
	The macroscopic density and velocity are obtained from the moments of the particle distribution functions as
	\begin{subequations}
		\begin{align}
			&\rho(\mathbf{x}_i,t)=\sum_{\alpha=0}^{8}f_\alpha(\mathbf{x}_i,t),\\
			&\mathbf{u}(\mathbf{x}_i,t)=\frac{1}{\rho(\mathbf{x}_i,t)}\sum_{\alpha=0}^{8}\mathbf{c}_\alpha f_\alpha(\mathbf{x}_i,t),
		\end{align}
	\end{subequations}
	and the lattice speed of sound is given by $c_s=\frac{\sqrt{3}}{3}$.
	
	In the numerical solution procedure, \eqref{LBE} is separated into the collision and streaming processes. The collision step is a local operation given by
	\begin{equation}\label{collision_step}
		f^*_\alpha(\mathbf{x}_i,t)=f_\alpha(\mathbf{x}_i,t)-\frac{1}{\tau}[f_\alpha(\mathbf{x}_i,t)-f_\alpha^{\rm eq}(\mathbf{x}_i,t)],
	\end{equation} 
	followed by the streaming step,
	\begin{equation}\label{streaming_step}
		f_\alpha(\mathbf{x}_i+\mathbf{c}_\alpha\Delta t,t+\Delta t)=f^*_\alpha(\mathbf{x}_i,t),
	\end{equation}
	where $f^*_\alpha(\mathbf{x}_i,t)$ represents the post-collision distribution.
	According to~\cite{Junk2005}, the kinematic viscosity is related to the relaxation time through
	\begin{equation}\label{viscosity_relation}
		\nu=\left(\tau-\frac{1}{2}\right)c_s^2.
	\end{equation}
	
	\subsubsection{Regularized Large-eddy simulation model}
	For unsteady flows containing vortical structures, it is generally impractical to resolve all spatial scales in numerical simulations. In the LES framework, flow structures larger than the filter width are directly resolved, while the effects of unresolved subgrid-scale motions are modeled through an eddy-viscosity closure. Following the Smagorinsky model~\cite{Smagorinsky1963}, the subgrid-scale eddy viscosity $\nu_t$ is given by
	\begin{equation}\label{nu_t}
		\nu_t=(C_s\ell_f)^2\|S\|,
	\end{equation}
	where $C_s$ is the Smagorinsky constant set to $0.1$, $\ell_f$ is the filter width, which is taken as the lattice spacing $\Delta x$ in the present LBM framework, and $\|S\|=\sqrt{2S_{pq}S_{pq}}$ denotes the magnitude of the strain-rate tensor, where repeated indices imply summation. Accordingly, an effective viscosity $\nu_{\rm eff}$ is defined as
	\begin{equation}\label{nu_eff}
		\nu_{\rm eff}=\nu+\nu_t.
	\end{equation}
	Using the viscosity-relaxation-time relation in \eqref{viscosity_relation}, the corresponding effective relaxation time is obtained as
	\begin{equation}\label{effective_viscosity_relation}
		\tau_{\rm eff}=\frac{\nu_{\rm eff}}{c_s^2}+\frac{1}{2}.
	\end{equation}
	Furthermore, according to the Chapman--Enskog analysis of the lattice Boltzmann equation~\cite{Lallemand2000,Hou1996}, the strain-rate tensor can be directly evaluated from the non-equilibrium moments of the distribution functions as
	\begin{equation}\label{Lallemand_equation}
		S_{pq}=-\frac{1}{2\rho c_s^2\tau_{\rm eff}}
		\sum_{\alpha=0}^{8}
		c_{\alpha p}c_{\alpha q}
		(f_\alpha-f_\alpha^{\rm eq}),
	\end{equation}
	where $c_{\alpha p}$ and $c_{\alpha q}$ denote the $p$-th and $q$-th components of the discrete velocity vector $\mathbf{c}_\alpha$, respectively, with $p,q\in\{1,2\}$ for the two-dimensional formulation.
	
	For convenience, the non-equilibrium second-order moment tensor $\Pi$ is introduced as
	\begin{equation}\label{Pi}
		\Pi_{pq}=
		\sum_{\alpha=0}^{8}
		c_{\alpha p}c_{\alpha q}
		\left(
		f_\alpha-f_\alpha^{\rm eq}
		\right),
	\end{equation}
	combining \eqref{nu_t}, \eqref{nu_eff}, \eqref{effective_viscosity_relation}, \eqref{Lallemand_equation}, and \eqref{Pi}, one obtains
	\begin{equation}\label{tau_eff}
		\tau_{\rm eff}=\frac12\left(\tau+\sqrt{\tau^2+\frac{2(C_s\ell_f)^2}{\rho c_s^4}\|\Pi\|}\right).
	\end{equation}
	
	Substituting the effective relaxation time \eqref{tau_eff} into the collision step \eqref{collision_step}, the collision operator can be rewritten as
	\begin{equation}\label{rewritten_collision}
		f_\alpha^*(\mathbf{x}_i,t)=f_\alpha^{\rm eq}(\mathbf{x}_i,t)+\left(1-\frac{1}{\tau_{\rm eff}}\right)f_\alpha^{\rm neq}(\mathbf{x}_i,t),
	\end{equation}
	where $f_\alpha^{\rm neq}(\mathbf{x}_i,t)=[f_\alpha(\mathbf{x}_i,t)-f_\alpha^{\rm eq}(\mathbf{x}_i,t)]$ denotes the non-equilibrium distribution function. To improve numerical stability, the regularization procedure of Latt and Chopard~\cite{Latt2006} is further applied. Only the first-order non-equilibrium contribution is retained, and the regularized collision step is expressed as
	\begin{equation}\label{modified_collision}
		f_\alpha^*(\mathbf{x}_i,t)=f_\alpha^{\rm eq}(\mathbf{x}_i,t)+\left(1-\frac{1}{\tau_{\rm eff}}\right)\frac{w_\alpha}{2c_s^4}
		\sum_{p,q}(c_{\alpha p}c_{\alpha q}-c_s^2\delta_{pq})\Pi_{pq},
	\end{equation}
	where $\delta_{pq}$ denotes the Kronecker delta. Therefore, the LES-regularized lattice Boltzmann equation used in the present work is given by
	\begin{equation}\label{modified_LES}
		f_\alpha(\mathbf{x}_i+\mathbf{c}_\alpha\Delta t,t+\Delta t)=f_\alpha^{\rm eq}(\mathbf{x}_i,t)+\left(1-\frac{1}{\tau_{\rm eff}}\right)\frac{w_\alpha}{2c_s^4}
		\sum_{p,q}(c_{\alpha p}c_{\alpha q}-c_s^2\delta_{pq})\Pi_{pq},
	\end{equation}
	which is implemented through the regularized collision step \eqref{modified_collision} followed by the streaming step \eqref{streaming_step}.
	
	\subsubsection{Initial and boundary conditions}
	Since the lattice Boltzmann equation is advanced in time, an initial condition must be prescribed for the particle distribution functions. In the present work, the distribution functions are initialized by the local equilibrium distribution corresponding to the given initial values of $\rho(\mathbf{x}_i,0), \mathbf{u}(\mathbf{x}_i,0)$, i.e.,
	\begin{equation}
		f_\alpha(\mathbf{x}_i,0)=f_\alpha^{\rm eq}(\rho(\mathbf{x}_i,0),\mathbf{u}(\mathbf{x}_i,0)).
	\end{equation}
	
	Let $\mathbf{n}$ denote the outward unit normal vector of the computational domain. During the streaming process, the distribution functions associated with lattice directions satisfying $\mathbf{c}_\alpha\cdot\mathbf{n}<0$ are unknown at the boundary nodes as illustrated in Fig.~\ref{illustration}. For no-slip wall boundaries, the bounce-back scheme is adopted to enforce the zero-velocity condition $\mathbf{u}=\mathbf{0}$. The unknown distribution functions are obtained through the bounce-back relation, namely
	\begin{equation}\label{wall}
		f_4=f_2,\quad f_7=f_5,\quad f_8=f_6.
	\end{equation}
	At the inlet and outlet boundaries, the Zou--He boundary condition~\cite{Zou1997} is adopted. This approach reconstructs the unknown distribution functions from the prescribed macroscopic variables based on the conservation of mass and momentum. Specifically, at the inlet
	\begin{equation}\label{inlet}
		\begin{aligned}
			f_1 &= f_3+\frac{2}{3}\rho u_0,\qquad\qquad\qquad\quad
			f_5 = f_7+\frac{1}{2}(f_4-f_2)
			+\frac{1}{6}\rho u_0,\\
			f_8 &= f_6+\frac{1}{2}(f_2-f_4)
			+\frac{1}{6}\rho u_0,\quad
			\rho = \frac{f_0+f_2+f_4+2(f_3+f_6+f_7)}{1-u_0},
		\end{aligned}
	\end{equation}
	where $u_0$ is the prescribed velocity in the $x$-direction at the inlet. For the outlet, the density is prescribed as $\rho_0$, and the unknown distribution functions are reconstructed as
	\begin{equation}\label{outlet}
		\begin{aligned}
			f_3 &= f_1-\frac{2}{3}\rho_0 u,\qquad\qquad\qquad\quad
			f_6 = f_8+\frac{1}{2}(f_4-f_2)
			-\frac{1}{6}\rho_0 u,\\
			f_7 &= f_5+\frac{1}{2}(f_2-f_4)
			-\frac{1}{6}\rho_0 u,\quad
			u = \frac{f_0+f_2+f_4+2(f_1+f_5+f_8)}{\rho_0}-1.
		\end{aligned}
	\end{equation}
	
	It should be noted that the configuration shown in Fig.~\ref{illustration} is presented only for illustration purposes. The specific flow configuration, boundary conditions, and solid geometry may vary depending on the problem under consideration. And the inlet velocity and outlet density are fixed at $u_0=0.03$ and $\rho_0=1.0$, respectively, throughout the present work.
	
	\subsection{Fluid-solid modeling}
	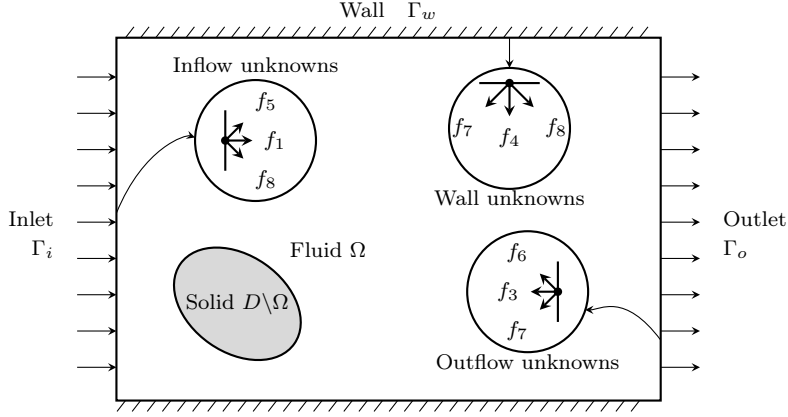
\begin{figure}[htbp]
		\centering
		\begin{tikzpicture}[scale=0.8,>=stealth]
			
			\draw[thick] (0,0) rectangle (9,6);
			
			\foreach \x in {0.15,0.4,...,8.85}{
				\draw[thin] (\x,6) -- ++(0.14,0.18);
				\draw[thin] (\x,0) -- ++(-0.14,-0.18);
			}
			\node at (4.5,6.45) {\scriptsize Wall\quad $\Gamma_w$};
			\node at (3.5,2.5){\scriptsize Fluid $\Omega$};
			
			\foreach \y in {0.55,1.15,...,5.45}{
				\draw[->,thin] (-0.65,\y) -- (0,\y);
				\draw[->,thin] (9,\y) -- (9.65,\y);
			}
			\node[left] at (-0.85,3) {\scriptsize Inlet};
			\node[left] at (-0.85,2.5) {\scriptsize $\Gamma_{i}$};
			
			\node[right] at (9.85,3) {\scriptsize Outlet};
			\node[right] at (9.85,2.5) {\scriptsize $\Gamma_{o}$};
			\draw[-,thick](1.8,3.8) -- (1.8,4.8);
			\draw[thick] (2.3,4.3) circle (1.0);
			\fill (1.8,4.3) circle (2pt);
			\draw[->,thick] (1.8,4.3) -- ++(0.45,0) node[right,font=\scriptsize] {$f_1$};
			\draw[->,thick] (1.8,4.3) -- ++(0.3,0.3) node[above right,font=\scriptsize] {$f_5$};
			\draw[->,thick] (1.8,4.3) -- ++(0.3,-0.3) node[below right,font=\scriptsize] {$f_8$};
			\node[above,font=\scriptsize] at (2.3,5.25) {Inflow unknowns};
			
			\draw[-,thick](7.3,1.3) -- (7.3,2.3);
			\draw[thick] (6.8,1.8) circle (1.0);
			\fill (7.3,1.8) circle (2pt);
			\draw[->,thick] (7.3,1.8) -- ++(-0.45,0)
			node[left,font=\scriptsize] {$f_3$};
			\draw[->,thick] (7.3,1.8) -- ++(-0.3,0.3)
			node[above left,font=\scriptsize] {$f_6$};
			\draw[->,thick] (7.3,1.8) -- ++(-0.3,-0.3)
			node[below left,font=\scriptsize] {$f_7$};
			\node[below,font=\scriptsize] at (6.8,0.95) {Outflow unknowns};
			
			\draw[->,thin] (0,3.1) .. controls (0.45,4.2) and (1.15,4.4) .. (1.3,4.35);
			\draw[->,thin] (9,1.0) .. controls (8.4,1.8) and (8.0,1.5) .. (7.75,1.5);
			
			\draw[-,thick] (6.0,5.25) -- (7.0,5.25);
			\draw[thick] (6.5,4.5) circle (1.0);
			\fill (6.5,5.25) circle (2pt);
			\draw[->,thick] (6.5,5.25) -- ++(0,-0.55)
			node[below,font=\scriptsize] {$f_4$};
			\draw[->,thick] (6.5,5.25) -- ++(-0.4,-0.4)
			node[below left,font=\scriptsize] {$f_7$};
			\draw[->,thick] (6.5,5.25) -- ++(0.4,-0.4)
			node[below right,font=\scriptsize] {$f_8$};
			\node[above,font=\scriptsize] at (6.5,3.05) {Wall unknowns};
			\draw[->,thin] (6.5,6.0) -- (6.5,5.5);
			
			\begin{scope}[shift={(2.0,1.6)}, rotate=-35]
				\draw[thick, fill=gray!30]
				(0,0) ellipse (1.15 and 0.80);
				\node at (0,0) {\scriptsize Solid\ $D\backslash\Omega$};
			\end{scope}
			
		\end{tikzpicture}
		\caption{Schematic illustration of the flow configuration.}
		\label{illustration}
	\end{figure}
	In topology optimization, a continuous design variable field $\phi$ is introduced to represent the material distribution within the computational domain. As illustrated in Fig.~\ref{illustration}, the computational domain $D$ consists of a fluid region $\Omega$ and a solid region $D\backslash\Omega$, and the design variable $\phi(\mathbf{x})$ is defined as
	\begin{equation}
		\phi(\mathbf{x}) =
		\begin{cases}
			1, & \text{if}\ \mathbf{x}\in\Omega,\\
			0, & \text{if}\ \mathbf{x}\in D\backslash\Omega.
		\end{cases}
	\end{equation}
	
	During the optimization process, the intermediate region between fluid and solid must be appropriately modeled. To this end, various interpolation schemes have been proposed, among which the Brinkman penalization approach is one of the most commonly used techniques~\cite{Borrvall2003}. To ensure that the introduced fluid-solid resistance preserves the first-order non-equilibrium information retained in \eqref{modified_LES}, the partial bounce-back formulation proposed by Zhu and Ma~\cite{Zhu2013} is adopted as
	\begin{equation}\label{external_force}
		F_\alpha(\mathbf{x}_i,t)=g(\phi)[f_{\bar{\alpha}}^{*}(\mathbf{x}_i,t)-f_\alpha^{*}(\mathbf{x}_i,t)],
	\end{equation}
	where $f_{\bar{\alpha}}^{*}$ is the post-collision value in the direction opposite $f_\alpha^{*}$, and $g(\phi)$ is a convex interpolation function adopted from N{\o}rgaard et al.~\cite{Norgaard2016},
	\begin{equation*}
		g(\phi)=0.5\bigg(1-\phi\left(\frac{1+\gamma}{\phi+\gamma}\right)\bigg),
	\end{equation*}
	where $\gamma$ is an adjustable parameter.
	\section{Representative-window adjoint framework}\label{Statistical_adjoint_framework}
	This section presents the construction of the proposed adjoint framework. A consecutive-window convergence criterion is first introduced to identify a representative time window within the fully developed flow stage. The objective function and discrete adjoint analysis are then formulated based on the identified time window, followed by the overall topology optimization procedure.
	\subsection{Consecutive-window convergence criterion}

	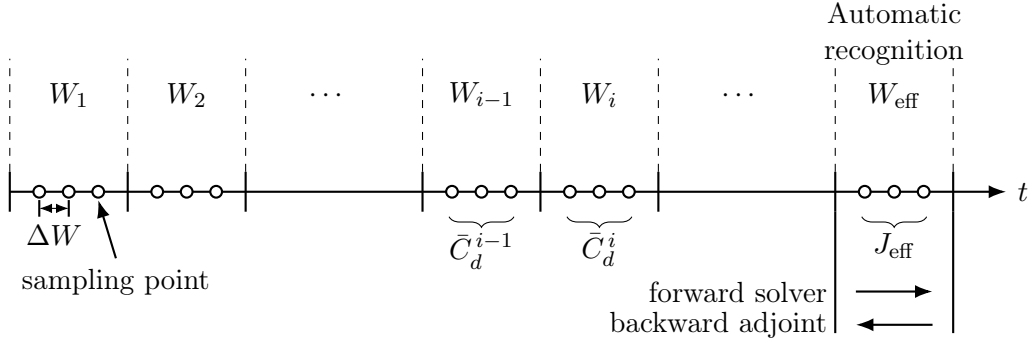
\begin{figure}[htbp]
		\centering
		\begin{tikzpicture}[
			x=0.78cm,
			y=2.2cm,
			>=Latex,
			every node/.style={font=\small}
			]
			
			\def\xzero{0}
			\def\xone{2.0}
			\def\xtwo{4.0}
			
			\def\xellipsisA{5.4}
			
			\def\ximinus{7.0}
			\def\xi{9.0}
			\def\xiright{11.0}
			
			\def\xellipsisB{12.4}
			
			\def\xeffL{14.0}
			\def\xeffR{16.0}
			
			\def\xend{16.9}
			
			\draw[->, thick]
			(0,0) -- (\xend,0)
			node[right] {$t$};
			
			\foreach \x in {
				\xzero,\xone,\xtwo,
				\ximinus,\xi,\xiright,
				\xeffL,\xeffR
			}{
				\draw[thick] (\x,-0.12) -- (\x,0.12);
				\draw[dashed] (\x,0.12) -- (\x,0.8);
			}
			
			\draw[thick] (\xeffL,-0.12) -- (\xeffL,-0.85);
			\draw[thick] (\xeffR,-0.12) -- (\xeffR,-0.85);
			\node at ({(\xzero+\xone)/2},0.58) {$W_1$};
			\node at ({(\xone+\xtwo)/2},0.58) {$W_2$};
			
			\node at ({(\ximinus+\xi)/2},0.58) {$W_{i-1}$};
			\node at ({(\xi+\xiright)/2},0.58) {$W_i$};
			
			\node at ({(\xeffL+\xeffR)/2},0.58) {$W_{\mathrm{eff}}$};
			\node[align=center] at ({(\xeffL+\xeffR)/2},0.95)
			{\(\text{Automatic}\)\\\(\text{recognition}\)};
			
			\node at (\xellipsisA,0.58) {$\cdots$};
			\node at (\xellipsisB,0.58) {$\cdots$};
			
			\foreach \x in {0.5,1.0,1.5}{
				\draw[fill=white, thick] (\x,0) circle (2.2pt);
			}
			
			\foreach \x in {2.5,3.0,3.5}{
				\draw[fill=white, thick] (\x,0) circle (2.2pt);
			}
			
			\foreach \x in {7.5,8.0,8.5}{
				\draw[fill=white, thick] (\x,0) circle (2.2pt);
			}
			
			\foreach \x in {9.5,10.0,10.5}{
				\draw[fill=white, thick] (\x,0) circle (2.2pt);
			}
			
			\foreach \x in {14.5,15.0,15.5}{
				\draw[fill=white, thick] (\x,0) circle (2.2pt);
			}
			
			\draw[->, thick](1.8,-0.4) -- (1.5,-0.05);
			
			\node[below]at (1.8,-0.4){sampling point};
			
			\draw[decorate,decoration={brace,mirror,amplitude=4pt}]
			(7.45,-0.15) -- (8.55,-0.15)
			node[midway,below=2pt] {$\bar{C}_d^{\,i-1}$};
			
			\draw[decorate,decoration={brace,mirror,amplitude=4pt}]
			(9.45,-0.15) -- (10.55,-0.15)
			node[midway,below=2pt] {$\bar{C}_d^{\,i}$};
			
			
			\draw[decorate,decoration={brace,mirror,amplitude=4pt}](14.45,-0.15) -- (15.55,-0.15)
			node[midway,below=2pt] {$J_{\mathrm{eff}}$};
			
			\draw[->, thick](14.35,-0.6) -- (15.65,-0.6);
			
			\node[left]at (14.05,-0.6){forward solver};
			
			\draw[<-, thick](14.35,-0.8) -- (15.65,-0.8);
			
			\node[left]at (14.05,-0.8){backward adjoint};
			
			\draw[<->,thin](0.5,-0.1) -- (1.0,-0.1)
			node[midway, below=2pt] {$\Delta W$};
			\draw[thick] (0.5,-0.15) -- (0.5,-0.05);
			\draw[thick] (1.0,-0.15) -- (1.0,-0.05);
			
		\end{tikzpicture}
		
		\caption{Schematic illustration of the consecutive time windows along the time evolution. }
		\label{statistical_windows}
	\end{figure}
	
	Motivated by the experimental observations discussed in Section~\ref{intro}, the temporal evolution of the flow is divided into a sequence of time windows. Let \(W_i\) denote the \(i\)-th time window. To determine whether the flow has reached the fully developed stage, the averaged value of a monitored quantity is evaluated and compared between two consecutive windows, \(W_{i-1}\) and \(W_i\).
	
	The mean drag coefficient is selected as the monitored averaged quantity in the present work. For a time window $W_i$, it is defined as
	\begin{equation}\label{statistical_force}
		\bar{C_d}^{i}=\frac{1}{|W_i|}\sum_{t \in W_i}C_d(t),
	\end{equation}
	where $|W_i|$ denotes the number of time steps within $W_i$, and $C_d(t)$ is the instantaneous drag coefficient at time $t$, given by
	\begin{equation}
		C_d(t)=\frac{\int_D F(\mathbf{x}_i,t) {\rm dx}}{\frac{1}{2}\rho_0 u_0^2 L_{\rm ref}},
	\end{equation}
	where $F(\mathbf{x}_i,t)$ denotes the first-order resistance force introduced in Eq.~\eqref{external_force}, and $L_{\rm ref}$ is a prescribed reference length. Based on the window-averaged drag coefficient, the consecutive-window convergence criterion is defined as
	\begin{equation}\label{statistical_criterion}
		\left|\frac{\bar{C_d}^{i}-\bar{C_d}^{i-1}}{\bar{C_d}^{i}}\right|<\varepsilon,
	\end{equation}
	where $\varepsilon$ denotes a prescribed tolerance. The overall forward procedure is schematically illustrated in Fig.~\ref{statistical_windows}. Once the criterion in Eq.~\eqref{statistical_criterion} is satisfied, the current time window $W_i$ is identified as the representative time window and is denoted by $W_{\rm eff}$.

	\subsection{Window-based discrete adjoint and sensitivity analysis}
	The resulting representative time window $W_{\rm eff}$ characterizes
	the fully developed flow stage associated with the current structure.
	The objective function is then evaluated over this window, consistent
	with the theoretical stability analysis of the cylinder wake performed
	after the initial transient has decayed~\cite{Dusek1994}.
	
	Accordingly, the objective function is defined as the following time-averaged quantity:
	\begin{equation}\label{obj}
		J_{\rm eff}(\mathbf{x};\phi)=\frac{1}{|W_{\rm eff}|}\sum_{t \in W_{\rm eff}} j(\mathbf{x},t;\phi),
	\end{equation}
	where $j(\mathbf{x},t;\phi)$ denotes the instantaneous objective function
	at spatial location $\mathbf{x}$ and time $t$. Within each design iteration, the identified \(W_{\mathrm{eff}}\) is treated as fixed when evaluating the objective and its adjoint sensitivity, which therefore corresponds to the current window-restricted subproblem. After the design update, \(W_{\mathrm{eff}}\) is re-identified from the new forward solution.
	
	To efficiently evaluate its sensitivity, the discrete adjoint method is employed. For notational simplicity, let superscripts $(i,n)$ denote the lattice node and time level, respectively. Accordingly, $f_\alpha(\mathbf{x}_i,t_n)$, $F_\alpha(\mathbf{x}_i,t_n)$ and $\lambda_\alpha(\mathbf{x}_i,t_n)$ are written as $f_\alpha^{i,n}$, $F_\alpha^{i,n}$ and $\lambda_\alpha^{i,n}$, respectively. With $\phi_i$ denoting the design variable at lattice node $i$, the discrete Lagrangian over $W_{\rm eff}$ is defined as
	\begin{equation}
		\begin{aligned}
			&\mathcal{L}(f_\alpha^{i,n},\lambda_\alpha^{i,n};\phi_i)\\
			&=J_{\rm eff}+\sum_{t_n\in W_{\rm eff}}\sum_{i,\alpha}{\lambda_\alpha^{i,n}\bigg(f_\alpha^{i+\alpha,n+1}-(f_\alpha^{\rm eq})^{i,n}-\left(1-\frac{1}{\tau_{\rm eff}}\right)(f_\alpha^{\rm neq})^{i,n}-F_\alpha^{i,n}\bigg)}.\\
		\end{aligned}
	\end{equation}
	Taking the partial derivative of the discrete Lagrangian with respect to the state variable $f_\alpha^{i,n}$ and imposing the stationarity condition
	\begin{equation}
		\frac{\partial \mathcal{L}}{\partial f_\alpha^{i,n}}=0,
	\end{equation}
	the resulting adjoint equation is
	\begin{equation}\label{adjoint_equation}
		\begin{aligned}
			&\lambda_\alpha^{i-\alpha,n-1}=-\frac{\partial J_{\rm eff}}{\partial f_\alpha^{i,n}}+\sum_{\beta=0}^{8}\lambda_\beta^{i,n}\bigg[\frac{\partial (f_\beta^{\rm eq})^{i,n}}{\partial f_\alpha^{i,n}}+\left(1-\frac{1}{\tau_{\rm eff}}\right)\frac{\partial (f_\beta^{\rm neq})^{i,n}}{\partial f_\alpha^{i,n}}\\
			&+\frac{1}{\tau_{\rm eff}^2}(f_\beta^{\rm neq})^{i,n}\frac{\partial\tau_{\rm eff}}{\partial f_\alpha^{i,n}}+g(\phi)\bigg(\frac{\partial((f_{\bar{\beta}}^{\rm eq})^{i,n}-(f_\beta^{\rm eq})^{i,n})}{\partial f_\alpha^{i,n}}\\
			&+\frac{1}{\tau_{\rm eff}^2}\left((f_{\bar{\beta}}^{\rm neq})^{i,n}-(f_\beta^{\rm neq})^{i,n}\right)\frac{\partial\tau_{\rm eff}}{\partial f_\alpha^{i,n}}+\left(1-\frac{1}{\tau_{\rm eff}}\right)\frac{\partial((f_{\bar{\beta}}^{\rm neq})^{i,n}-(f_\beta^{\rm neq})^{i,n})}{\partial f_\alpha^{i,n}}\bigg)\bigg],
		\end{aligned}
	\end{equation}
	where,
	\begin{equation}
		\begin{aligned}
			&\frac{\partial(f_\beta^{\rm eq})^{i,n}}{\partial f_\alpha^{i,n}}=w_\beta\bigg[1+\frac{\mathbf{c_\beta}\cdot\mathbf{u}^{i,n}}{c_s^2}+\frac{(\mathbf{c}_\beta\cdot\mathbf{u}^{i,n})^2}{2c_s^4}-\frac{\mathbf{u}^{i,n}\cdot\mathbf{u}^{i,n}}{2c_s^2}\bigg]\\
			&+w_\beta\bigg[\frac{(\mathbf{c_\beta}-\mathbf{u}^{i,n})\cdot(\mathbf{c}_\alpha-\mathbf{u}^{i,n})}{c_s^2}+\frac{(\mathbf{c}_\beta\cdot\mathbf{u}^{i,n})(\mathbf{c}_\beta\cdot(\mathbf{c}_\alpha-\mathbf{u}^{i,n}))}{c_s^4}\bigg],\\
		\end{aligned}
	\end{equation}
	\begin{equation}
		\begin{aligned}
			&\frac{\partial \tau_{\rm eff}}{\partial f_\alpha^{i,n}}=\frac{(C_s\ell_f)^2}{2c_s^4\sqrt{\tau^2+\frac{2(C_s\ell_f)^2}{\rho^{i,n} c_s^4}\|\Pi^{i,n}\|}}\left(-\frac{\|\Pi^{i,n}\|}{(\rho^{i,n})^2}+\frac{1}{\rho^{i,n}}\frac{\partial \|\Pi^{i,n}\|}{\partial f_\alpha^{i,n}}\right),\\
		\end{aligned}
	\end{equation}
	\begin{equation}
		\begin{aligned}
			&\frac{\partial\Pi_{pq}^{i,n}}{\partial f_\alpha^{i,n}}=-\sum_{\beta=0}^8c_{\beta p}c_{\beta q}\frac{\partial (f_\beta^{\rm eq})^{i,n}}{\partial f_\alpha^{i,n}}+c_{\alpha p}c_{\alpha q},
		\end{aligned}
	\end{equation}
	\begin{equation}
		\begin{aligned}
			&\frac{\partial(f_\beta^{\rm neq})^{i,n}}{\partial f_\alpha^{i,n}}=\frac{w_\beta}{2c_s^4}\sum_{p,q}\left(c_{\beta p}c_{\beta q}-c_s^2\delta_{pq}\right)\frac{\partial \Pi_{pq}^{i,n}}{\partial f_\alpha^{i,n}}.\\
		\end{aligned}
	\end{equation}
	Furthermore, \eqref{adjoint_equation} can be decomposed into an adjoint collision step
	\begin{equation}\label{adjoint_collision}
		\begin{aligned}
			\lambda_\alpha^{i,n-1}&=-\frac{\partial J_{\rm eff}}{\partial f_\alpha^{i,n}}+\sum_{\beta=0}^{8}\lambda_\beta^{i,n}\bigg[\frac{\partial (f_\beta^{\rm eq})^{i,n}}{\partial f_\alpha^{i,n}}+\left(1-\frac{1}{\tau_{\rm eff}}\right)\frac{\partial (f_\beta^{\rm neq})^{i,n}}{\partial f_\alpha^{i,n}}\\
			&+\frac{1}{\tau_{\rm eff}^2}(f_\beta^{\rm neq})^{i,n}\frac{\partial\tau_{\rm eff}}{\partial f_\alpha^{i,n}}+g(\phi)\bigg(\frac{\partial((f_{\bar{\beta}}^{\rm eq})^{i,n}-(f_\beta^{\rm eq})^{i,n})}{\partial f_\alpha^{i,n}}\\
			&+\frac{1}{\tau_{\rm eff}^2}\left((f_{\bar{\beta}}^{\rm neq})^{i,n}-(f_\beta^{\rm neq})^{i,n}\right)\frac{\partial\tau_{\rm eff}}{\partial f_\alpha^{i,n}}+\left(1-\frac{1}{\tau_{\rm eff}}\right)\frac{\partial((f_{\bar{\beta}}^{\rm neq})^{i,n}-(f_\beta^{\rm neq})^{i,n})}{\partial f_\alpha^{i,n}}\bigg)\bigg],
		\end{aligned}
	\end{equation}
	and an adjoint streaming step
	\begin{equation}\label{adjoint_streaming}
		\lambda_\alpha^{i-\alpha,n-1}=\lambda_\alpha^{i,n-1}.
	\end{equation}
	
	The adjoint equations are solved backward in time from the end of the representative time window $W_{\rm eff}$ to its beginning as schematically illustrated in Fig.~\ref{statistical_windows}.
	
	Denote \(\widehat{\lambda}_\alpha\) to be the adjoint distribution functions after the boundary treatment. Following the adjoint boundary treatment of Yaji et al.~\cite{Yaji2014}, the boundary conditions are discretized consistently with the D2Q9 lattice and the corresponding forward boundary schemes. For the no-slip wall boundary, the adjoint bounce-back condition is given by
	\begin{equation}
		\widehat{\lambda}_2=\lambda_4,\qquad
		\widehat{\lambda}_5=\lambda_7,\qquad
		\widehat{\lambda}_6=\lambda_8.
	\end{equation}
	For conciseness, we introduce $\chi=\frac{u_0}{1-u_0}$. For the inlet boundary, the adjoint Zou--He boundary condition is given by
	\begin{equation}
		\begin{aligned}
			\widehat{\lambda}_0
			&=\lambda_0
			+\frac{2\chi}{3}\lambda_1
			+\frac{\chi}{6}\lambda_5
			+\frac{\chi}{6}\lambda_8,\\
			\widehat{\lambda}_2
			&=\lambda_2
			+\frac{2\chi}{3}\lambda_1
			+\left(-\frac{1}{2}+\frac{\chi}{6}\right)\lambda_5
			+\left(\frac{1}{2}+\frac{\chi}{6}\right)\lambda_8,\\
			\widehat{\lambda}_4
			&=\lambda_4
			+\frac{2\chi}{3}\lambda_1
			+\left(\frac{1}{2}+\frac{\chi}{6}\right)\lambda_5
			+\left(-\frac{1}{2}+\frac{\chi}{6}\right)\lambda_8,\\
			\widehat{\lambda}_3
			&=\lambda_3
			+\left(1+\frac{4\chi}{3}\right)\lambda_1
			+\frac{\chi}{3}\lambda_5
			+\frac{\chi}{3}\lambda_8,\\
			\widehat{\lambda}_6
			&=\lambda_6
			+\frac{4\chi}{3}\lambda_1
			+\frac{\chi}{3}\lambda_5
			+\left(1+\frac{\chi}{3}\right)\lambda_8,\\
			\widehat{\lambda}_7
			&=\lambda_7
			+\frac{4\chi}{3}\lambda_1
			+\left(1+\frac{\chi}{3}\right)\lambda_5
			+\frac{\chi}{3}\lambda_8,\\
			\widehat{\lambda}_1
			&=\widehat{\lambda}_5
			=\widehat{\lambda}_8=0.
		\end{aligned}
		\label{eq:adjoint_inlet_bc}
	\end{equation}
	
	For the outlet boundary, the adjoint Zou--He boundary condition is
	given by
	\begin{equation}
		\begin{aligned}
			\widehat{\lambda}_0
			&=\lambda_0
			-\frac{2}{3}\lambda_3
			-\frac{1}{6}\lambda_6
			-\frac{1}{6}\lambda_7,\\
			\widehat{\lambda}_1
			&=\lambda_1
			-\frac{1}{3}\lambda_3
			-\frac{1}{3}\lambda_6
			-\frac{1}{3}\lambda_7,\\
			\widehat{\lambda}_2
			&=\lambda_2
			-\frac{2}{3}\lambda_3
			-\frac{2}{3}\lambda_6
			+\frac{1}{3}\lambda_7,\\
			\widehat{\lambda}_4
			&=\lambda_4
			-\frac{2}{3}\lambda_3
			+\frac{1}{3}\lambda_6
			-\frac{2}{3}\lambda_7,\\
			\widehat{\lambda}_5
			&=\lambda_5
			-\frac{4}{3}\lambda_3
			-\frac{1}{3}\lambda_6
			+\frac{2}{3}\lambda_7,\\
			\widehat{\lambda}_8
			&=\lambda_8
			-\frac{4}{3}\lambda_3
			+\frac{2}{3}\lambda_6
			-\frac{1}{3}\lambda_7,\\
			\widehat{\lambda}_3
			&=\widehat{\lambda}_6
			=\widehat{\lambda}_7=0.
		\end{aligned}
		\label{eq:adjoint_outlet_bc}
	\end{equation}
	After solving the forward and adjoint systems, the sensitivity at the $i$-th lattice node is obtained by
	\begin{equation}\label{sensitivity}
		\begin{aligned}
			&\frac{\partial \mathcal{L}}{\partial \phi_i}=\frac{\partial J_{\rm eff}}{\partial \phi_i}-\sum_{t_n\in W_{\rm eff}}\sum_{\alpha}\frac{\gamma(1+\gamma)}{2(\phi_i+\gamma)^2}\lambda_\alpha^{i,n}\bigg(((f_\alpha^{\rm eq})^{i,n}-(f_{\bar{\alpha}}^{\rm eq})^{i,n})\\
			&+\left(1-\frac{1}{\tau_{\rm eff}}\right)((f_\alpha^{\rm neq})^{i,n}-(f_{\bar{\alpha}}^{\rm neq})^{i,n})\bigg).\\
		\end{aligned}
	\end{equation}
	\subsection{Topology optimization procedure}
	During the optimization process, the design variables are regularized using a density filter with radius \(r=4\Delta x\) and mapped to the physical material field \(\phi_i\) through the Heaviside projection~\cite{Wang2011}. The projection threshold is set to \(\eta=0.5\), while the projection parameter \(\beta\) is initialized at 2 and doubled every 10 optimization iterations until reaching its maximum value of 16. The volume fraction is defined as
	\begin{equation}
		V=\frac{1}{n_D}\sum_{i=1}^{n_D}(1-\phi_i),
	\end{equation}
	where \(n_D\) denotes the total number of lattice nodes in the design domain. The design variables are updated using the Method of Moving Asymptotes (MMA)~\cite{Svanberg1987} under the prescribed volume constraint. Before being supplied to the MMA optimizer, the sensitivities obtained from Eq.~\eqref{sensitivity} are transformed through the density filter and Heaviside projection using the standard chain-rule treatment in~\cite{Wang2011}. The optimization terminates when either \(N_{\max}\) is reached or the relative change in \(J_{\rm eff}\) falls below the prescribed tolerance \(\varepsilon_o\).
	
	To evaluate the window-averaged quantities, a temporal sampling strategy
	is adopted within each time window. Let $\Delta W$ denote the interval
	between two consecutive sampling instants, as illustrated in
	Fig.~\ref{statistical_windows}. The averaged quantities in
	Eqs.~\eqref{statistical_force} and~\eqref{obj} are evaluated using the
	weighted samples as
	\begin{equation}
		\overline{C}_d^i
		=
		\frac{\Delta W}{|W_i|}
		\sum_{t_k\in W_i} C_d(t_k),
		\qquad
		J_{\rm eff}(\boldsymbol{x};\phi)
		=
		\frac{\Delta W}{|W_{\rm eff}|}
		\sum_{t_k\in W_{\rm eff}}
		j(\boldsymbol{x},t_k;\phi),
		\label{weighted_temporal_sampling}
	\end{equation}
	where $t_k$ denotes the sampling instants. Consequently, the objective
	source term $-\partial J_{\rm eff}/\partial f_\alpha^{i,n}$ in
	Eq.~\eqref{adjoint_collision} is nonzero only at the sampling instants
	and incorporates the corresponding temporal weight. It vanishes at the
	intermediate time steps, whereas the adjoint equations are solved backward
	at every time step within $W_{\rm eff}$.
	
	In summary, the computational procedures of the regularized LES solver and the topology optimization framework are summarized in Algorithm~\ref{LES_agri} and Algorithm~\ref{TO_agri}, respectively.
	
	\begin{algorithm}[H]
		\caption{LES solver}
		\label{LES_agri}
		
		Initialization.
		
		\For{Time window $W_i$, $i=1,2,\cdots,W_{\rm max}$}{
			
			\For{each time step $t\in W_i$}{
				
				Collision:
				
				\Indp
				Compute the equilibrium distribution function
				$f_\alpha^{\rm eq}$ by Eq.~\eqref{feq}.
				
				Compute the non-equilibrium second-order moment tensor
				$\Pi$ by Eq.~\eqref{Pi}.
				
				Compute the effective relaxation time
				$\tau_{\rm eff}$ by Eq.~\eqref{tau_eff}.
				
				Perform the regularized collision by
				Eq.~\eqref{modified_collision}.
				
				Apply the fluid-solid treatment by
				Eq.~\eqref{external_force}.
				\Indm
				
				Streaming by Eq.~\eqref{streaming_step}.
				
				Apply boundary conditions by
				Eqs.~\eqref{wall}--\eqref{outlet}.
				
				\If{the current time step is a sampling instant}{
					
					Update the weighted averaged quantity
					$\overline{C}_d^i$ by
					Eq.~\eqref{weighted_temporal_sampling}.
				}
			}
			
			\If{convergence is confirmed according to
				Eq.~\eqref{statistical_criterion}}{
				
				Set $W_{\rm eff}=W_i$.
				
				\Return{$W_{\rm eff}$ and the corresponding
					distribution functions $f_\alpha$.}
			}
		}
	\end{algorithm}
	\begin{algorithm}[H]
		\caption{Optimization iteration}
		\label{TO_agri}
		
		Initialize the design variable $\phi_0$.
		
		\While{the iteration number $m<N_{\max}$}{
			
			Solve the forward problem using Algorithm~\ref{LES_agri} to obtain
			$W_{\rm eff}$ and the corresponding distribution functions $\{f_\alpha\}$.
			
			\For{each time step $t_n\in W_{\rm eff}$ in reverse chronological order}{
				
				Solve the adjoint equations using
				Eqs.~\eqref{adjoint_collision} and~\eqref{adjoint_streaming}.
			}
			
			Compute the design sensitivities based on Eq.~\eqref{sensitivity}.
			
			Apply the density filter.
			
			Apply the Heaviside projection.
			
			Update the design variables using MMA.
			
			\If{the relative change in $J_{\rm eff}$ is smaller than
				$\varepsilon_o$}{
				\textbf{break}.
			}
		}
		
		\Return{Optimized result.}
		
	\end{algorithm}
	\section{Numerical Results}\label{Numerical_Results}
	In this section, the backward-facing step (BFS) is first employed to validate the accuracy of the forward solver. Subsequently, the flow around a circular cylinder is considered to investigate the characteristics of the proposed approach in topology optimization. The wake flow recovery problem is then employed to demonstrate its capability to identify the fully developed stage of unsteady flows. Finally, a U-bend flow is considered to further evaluate its performance in more complex unsteady flows.
	
	In all numerical examples, the Reynolds number is defined as
	\begin{equation}
		{\rm Re}=\frac{UL}{\nu},
	\end{equation}
	where $U$ denotes the characteristic velocity and $L$ denotes the characteristic length associated with the corresponding flow configuration.
	
	\subsection{Backward-facing step}
	\begin{figure}[htbp]
		\centering
		\begin{tikzpicture}[scale=0.78,>=Stealth]
			
			\draw[thick] (-4.0,1.2) -- (6.2,1.2);
			\draw[thick] (-4.0,-0.3) -- (0,-0.3);
			\draw[thick] (0,-0.3) -- (0,-1.2);
			\draw[thick] (0,-1.2) -- (6.2,-1.2);
			
			\draw[<->] (-4.0,1.4) -- (0.0,1.4);
			\node[above] at (-2.0,1.4){$L_u$};
			\draw[<->] (0.0,1.4) -- (6.2,1.4);
			\node[above] at (3.1,1.4){$L_d$};
			\draw[->,thick] (-3.8,0.45) -- (-2.7,0.45);
			\draw[->,thick] (-3.8,0.75) -- (-2.7,0.75);
			\draw[->,thick] (-3.8,0.15) -- (-2.7,0.15);
			
			\draw[->,thick] (5.9,0.0) -- (6.7,0.0);
			\draw[->,thick] (5.9,0.45) -- (6.7,0.45);
			\draw[->,thick] (5.9,-0.45) -- (6.7,-0.45);
			
			
			\coordinate (O) at (0,-0.3);
			\draw[->] (O) -- ++(0.75,0);
			\draw[->] (O) -- ++(0,0.75);
			\node[right] at (0.75,-0.3) {$x$};
			\node[above] at (0.4,0.05) {$y$};
			
			\draw[thick,->]
			(0,-0.3)
			.. controls (0.7,-0.35) and (1.4,-0.65)
			.. (2.2,-1.2);
			\draw[thick,->]
			(0.25,-0.95)
			arc[start angle=210,end angle=-120,x radius=0.5,y radius=0.25];
			
			\draw[thick]
			(0.25,-0.95)
			arc[start angle=210,end angle=500,x radius=0.5,y radius=0.25];
			
			\draw[thick]
			(2.45,1.2)
			arc[start angle=180,end angle=360,
			x radius=0.75,
			y radius=0.48];
			
			\draw[thick,->]
			(2.9,0.92)
			arc[start angle=210,end angle=-120,
			x radius=0.35,
			y radius=0.11];
			
			\draw[thick]
			(2.9,0.92)
			arc[start angle=210,end angle=500,
			x radius=0.35,
			y radius=0.11];
			
			\draw[<->] (-0.35,-0.3) -- (-0.35,-1.2);
			\node[left] at (-0.35,-0.75) {$H$};
			
			\draw[<->] (-0.35,1.2) -- (-0.35,-0.3);
			\node[left] at (-0.35,0.45) {$H_u$};
			
			\draw[<->] (5.25,1.2) -- (5.25,-1.2);
			\node[left] at (5.25,0.45) {$H_d$};
			
			\draw[<->] (0,-1.55) -- (2.2,-1.55);
			\node[below] at (1.1,-1.55) {$x_r$};
			
			\node[above] at (5.45,-2.0) {Wall};
			
		\end{tikzpicture}
		\caption{Schematic diagram of the backward-facing step flow.}
		\label{BFS}
	\end{figure}
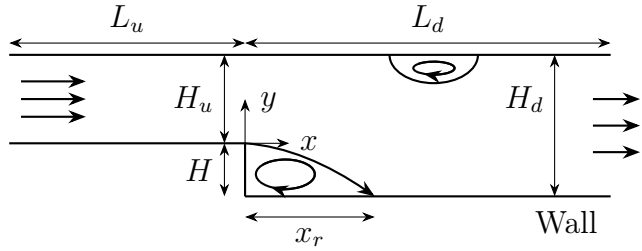
	
	The BFS problem is first employed to validate the accuracy of the forward solver. As a classical benchmark for separated internal flows, this problem has been extensively investigated experimentally and theoretically by Armaly et al.~\cite{Armaly1983}. In their numerical study, the flow was solved using a two-dimensional steady finite-difference solver for the incompressible Navier--Stokes equations, while the experimental measurements were used as the primary reference for validation. Experimental results show that the flow remains essentially two-dimensional for Reynolds numbers up to approximately ${\rm Re}=400$, while noticeable three-dimensional effects appear at higher Reynolds numbers. Therefore, the present two-dimensional solver is validated for $100 \leq {\rm Re} \leq 400$, where direct comparisons with the reported experimental data are performed.
	
	The computational configuration is shown in Fig.~\ref{BFS}. Here, $H$ denotes the step height, $H_u$ and $H_d$ are the upstream and downstream channel heights, respectively, $L_u$ and $L_d$ denote the upstream and downstream lengths, respectively, and $x_r$ represents the reattachment length measured from the step corner. The origin of the coordinate system is located at the step corner.
	\begin{table}[htbp]
		\centering
		\begin{tabular}{lc}
			\hline
			Parameter & Value \\
			\hline
			$H_u$ & 30 \\
			$H$ & 30 \\
			$H_d$ & 60 \\
			$L_u$ & 300 \\
			$L_d$ & 1200 \\
			$\varepsilon$ & $1.0\times10^{-3}$ \\
			$|W_i|$ & 5000 \\
			$\Delta W$ & 100\\
			\hline
		\end{tabular}
		\caption{Computational parameters for the backward-facing step problem.}
		\label{BFS_parameters}
	\end{table}
	
\begin{figure}[htbp]
		\centering
		\includegraphics[width=0.5\linewidth]{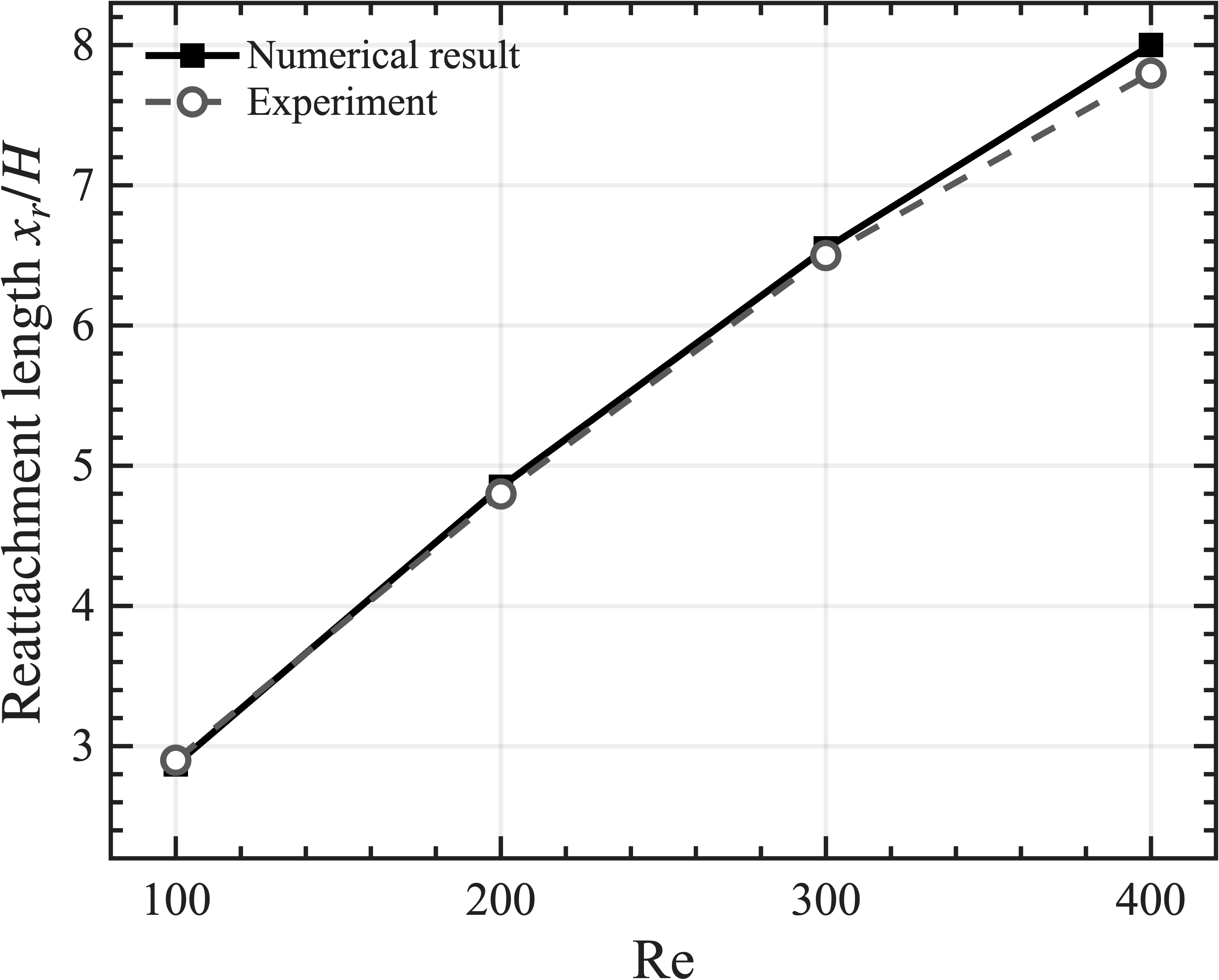}
		\caption{Comparison of the normalized reattachment length with the experimental data (originally Fig. 13(a) in Armaly et al.~\cite{Armaly1983}).}
		\label{reattachmentlength}
	\end{figure}

	The computational parameters are summarized in Table~\ref{BFS_parameters}. The predicted reattachment lengths at different Reynolds numbers are compared with the corresponding experimental results, as shown in Fig.~\ref{reattachmentlength}. Excellent agreement is observed between the numerical results and the experimental data.
	
\begin{figure}[htbp]
		\centering
		
		\begin{subfigure}{0.48\linewidth}
			\centering
			\includegraphics[width=\linewidth]{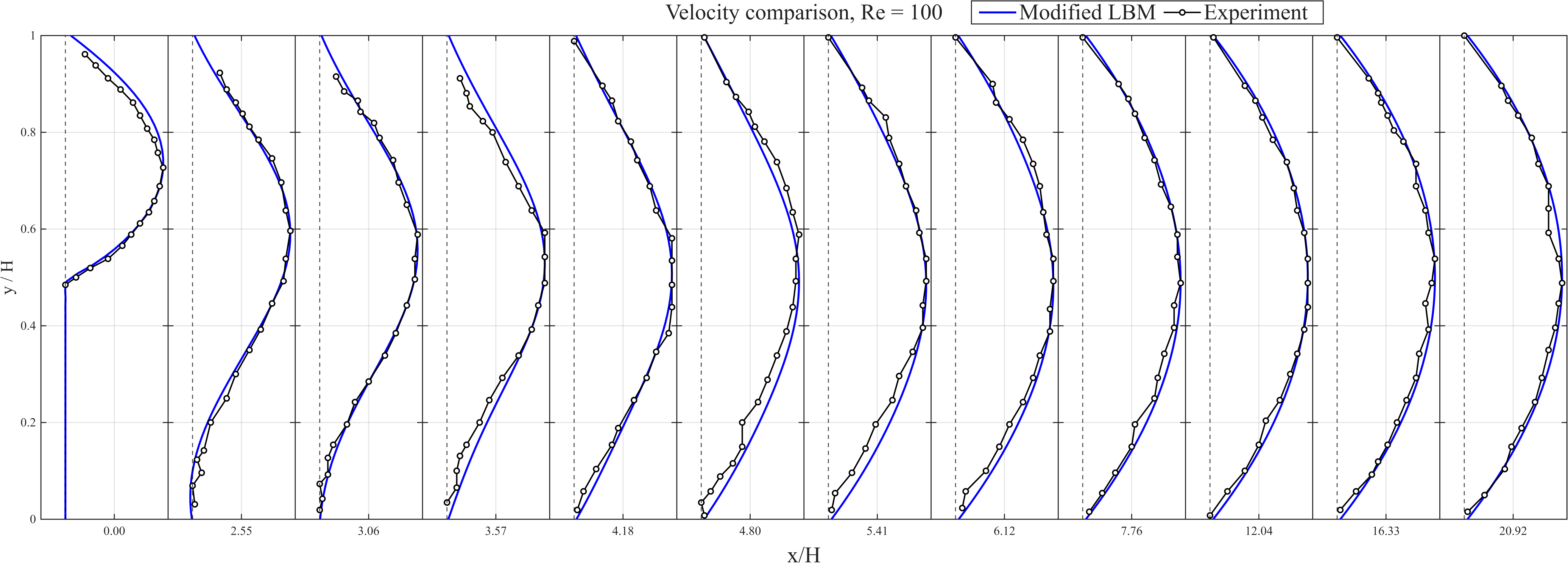}
			\caption{Streamwise velocity profiles, ${\rm Re}=100$.}
			\label{armalyexpvslbmwideprofilesre100}
		\end{subfigure}
		\hfill
		\begin{subfigure}{0.48\linewidth}
			\centering
			\includegraphics[width=\linewidth]{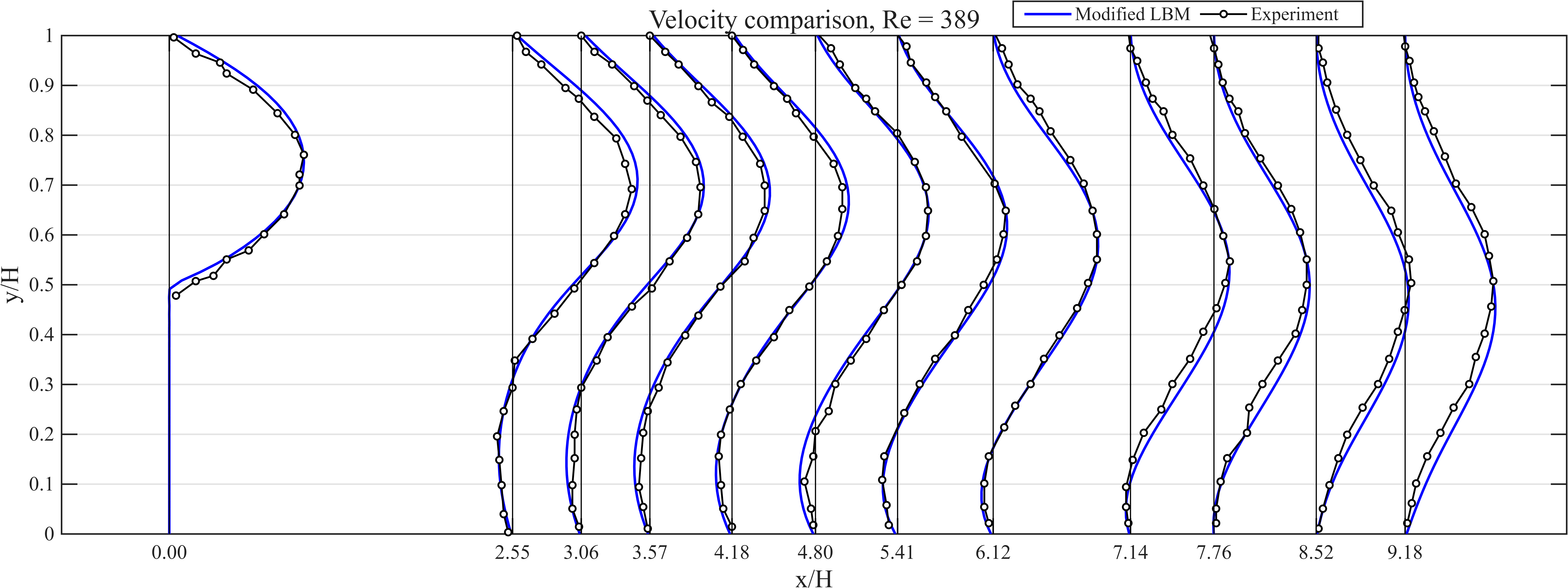}
			\caption{Streamwise velocity profiles, ${\rm Re}=389$.}
			\label{armalyexpvslbmwideprofilesre389}
		\end{subfigure}
		
		\caption{Comparison of the streamwise velocity profiles with the experimental data (originally Fig. 15 and 16 in Armaly et al.~\cite{Armaly1983}).}
		\label{velocityprofile}
	\end{figure}

\begin{figure}[htbp]
		\centering
		\begin{subfigure}{0.45\linewidth}
			\centering
			\includegraphics[width=\linewidth]{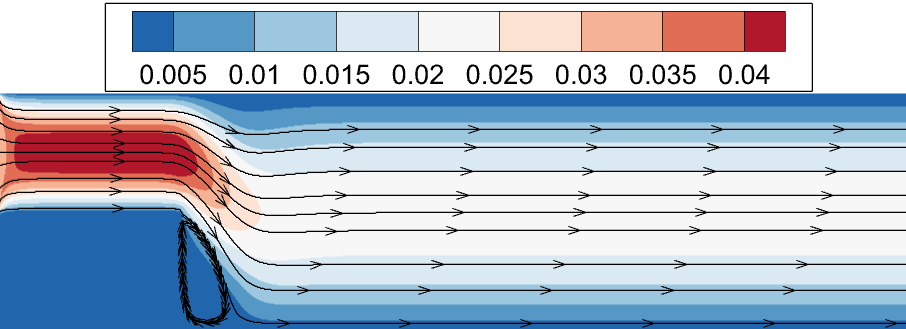}
			\caption{Velocity field, ${\rm Re}=100$.}
			\label{velocityfieldre100}
		\end{subfigure}
		\hfill
		\begin{subfigure}{0.45\linewidth}
			\centering
			\includegraphics[width=\linewidth]{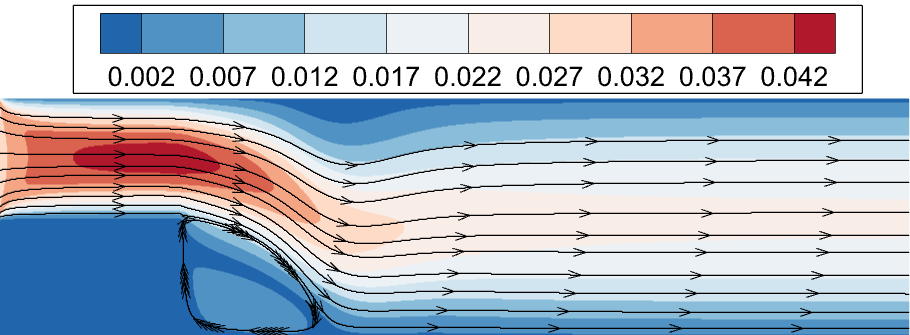}
			\caption{Velocity field, ${\rm Re}=389$.}
			\label{velocityfieldre389}
		\end{subfigure}
		
		\caption{Results of different velocity fields.}
		\label{velocityfield}
	\end{figure}

	To further assess the accuracy of the present solver, the velocity profiles at different downstream locations are compared with the experimental measurements for ${\rm Re}=100$ and ${\rm Re}=389$. Fig.~\ref{velocityprofile} shows the comparison of the streamwise velocity profiles at several cross-sections located downstream of the step, which shows good agreement with the experimental data. The flow fields obtained when the consecutive-window convergence criterion is satisfied are shown in Fig.~\ref{velocityfield}.
	
	
	Overall, the present solver accurately predicts both the global reattachment characteristics and the local velocity distributions of the BFS flow, thereby providing a reliable foundation for the subsequent topology optimization studies. It also indicates that the consecutive-window convergence criterion terminates the time-marching procedure only after the flow has reached a sufficiently developed stage, thereby providing a physically reasonable stopping condition for the present benchmark.
	
	\subsection{Flow around a circular cylinder}
	Following the validation of the forward solver, the flow around a circular cylinder is considered as a representative benchmark to validate the proposed framework. Owing to its rich flow dynamics, this problem has been extensively investigated experimentally, and the evolution of wake patterns over a wide range of Reynolds numbers has been systematically summarized by Williamson~\cite{Williamson1996}. The study showed that the wake undergoes a sequence of distinct flow regimes as the Reynolds number increases. 
	
	At very low Reynolds numbers, the flow remains attached to the cylinder without noticeable wake separation. As the Reynolds number increases, a pair of steady symmetric recirculation bubbles gradually develops behind the cylinder. When the Reynolds number exceeds the critical value of approximately ${\rm Re}=49$, a Hopf bifurcation occurs and the flow transitions to a laminar vortex-shedding regime characterized by periodic wake oscillations. Moreover, the wake remains essentially two-dimensional up to approximately ${\rm Re}=190$, beyond which significant three-dimensional instabilities develop.
	
	Accordingly, three Reynolds numbers, ${\rm Re}=3.2$, $12.8$, and $51.2$, are considered in the present study. These cases respectively represent the attached-flow regime, the steady symmetric-wake regime, and the vortex-shedding regime, thereby providing a systematic assessment of the proposed representative-window adjoint framework under different flow mechanisms. A $200\times100$ lattice is adopted, where a cylinder of radius $8$ is placed at the center of the computational domain.
	
\begin{figure}[htbp]
		\centering
		
		\begin{subfigure}{0.32\linewidth}
			\centering
			\includegraphics[width=\linewidth]{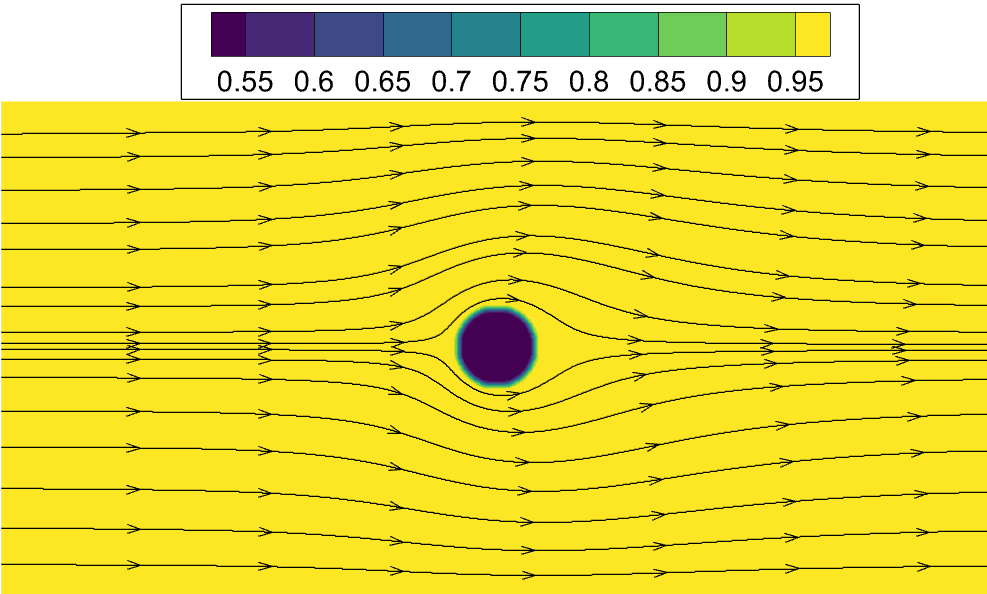}
			\caption{${\rm Re}=3.2$.}
			\label{initialre3}
		\end{subfigure}
		\hfill
		\begin{subfigure}{0.32\linewidth}
			\centering
			\includegraphics[width=\linewidth]{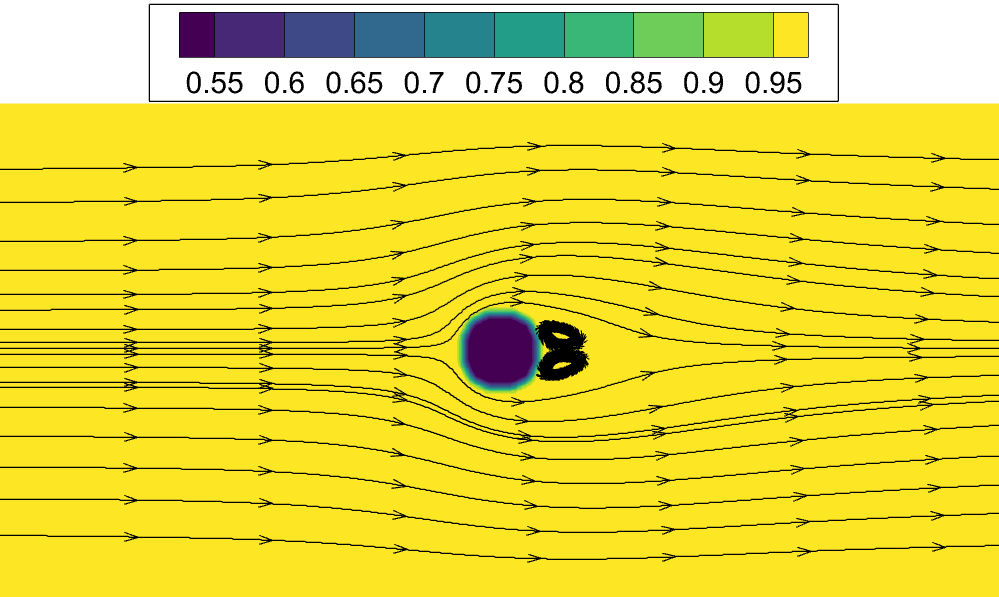}
			\caption{${\rm Re}=12.8$.}
			\label{initialre12}
		\end{subfigure}
		\hfill
		\begin{subfigure}{0.32\linewidth}
			\centering
			\includegraphics[width=\linewidth]{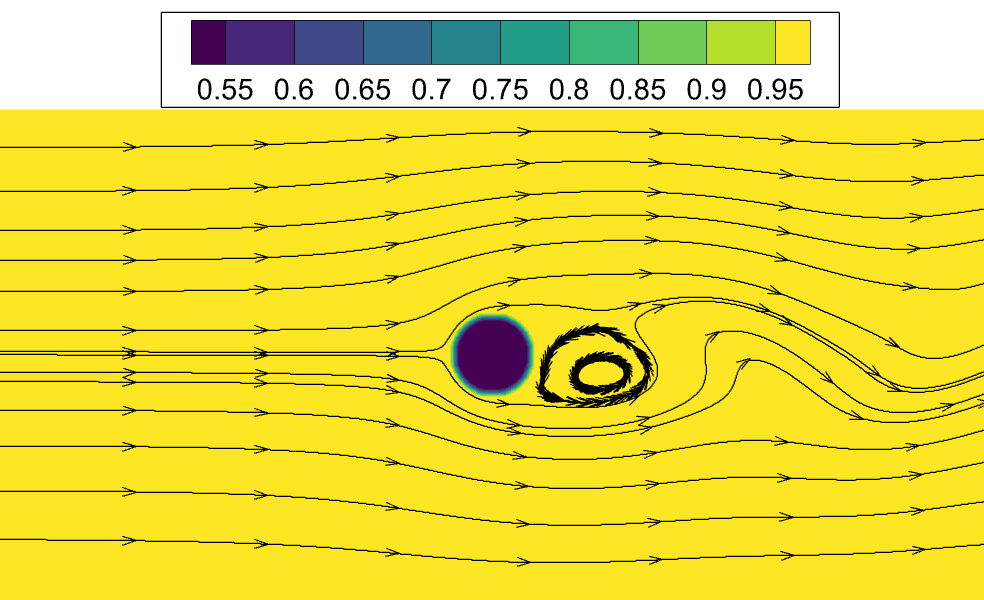}
			\caption{${\rm Re}=51.2$.}
			\label{initialre51}
		\end{subfigure}
		
		\caption{Flow fields for different Reynolds numbers in design domain.}
		\label{initial_flow_fields}
	\end{figure}

	The objective is to minimize averaged drag force over $W_{\rm eff}$, defined as
	\begin{equation}
		\begin{aligned}
			&\min_{\phi}F_x=\frac{\Delta W}{|W_{\rm eff}|}\sum_{t\in W_{\rm eff}}\sum_{i\in D}\sum_{\alpha=0}^8 c_{x,\alpha}F_\alpha(\mathbf{x}_i,t),\\
			&\text{s.t.}\quad 0.014 \leq V \leq 0.016 .
		\end{aligned}
	\end{equation}
	where $c_{x,\alpha}$ denotes the $x$-component of the discrete velocity vector $\mathbf{c}_\alpha$. The tolerance $\varepsilon$ is set to be $10^{-3}$ and the initial design variable inside the cylinder is set to $\phi=0.5$. The corresponding optimization problem is known to admit a streamlined (teardrop-shaped) result for drag reduction, providing a suitable benchmark for evaluating the proposed framework. The initial distribution and flow fields of different $\rm Re$ are shown in Fig.~\ref{initial_flow_fields}.
	
	\subsubsection{Single-point sensitivity validation}
\begin{figure}[htbp]
		\centering
		
		\begin{subfigure}{0.32\linewidth}
			\centering
			\includegraphics[width=\linewidth]{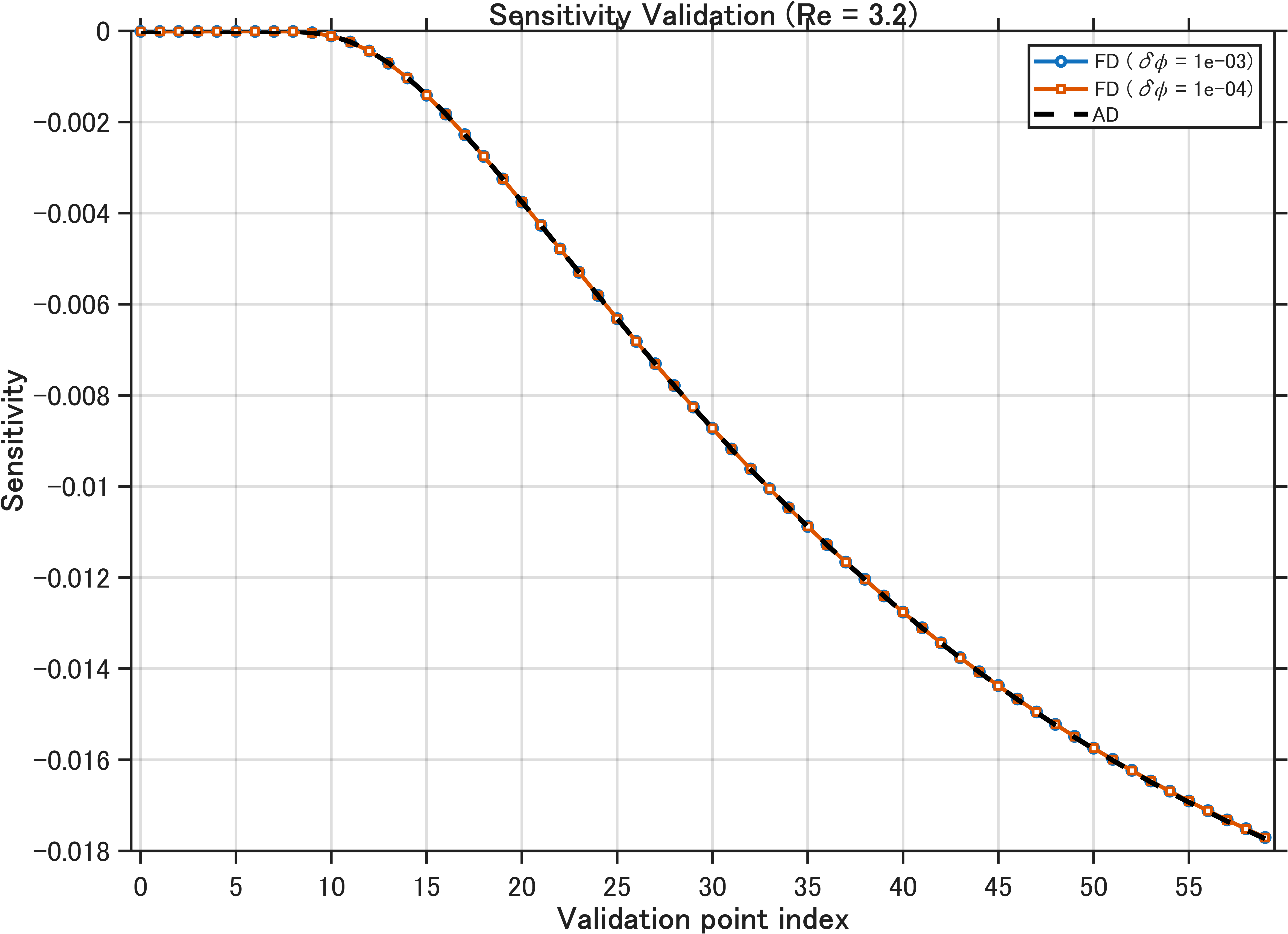}
			\caption{${\rm Re}=3.2$.}
			\label{sensitivityvalidationdifferentepsre3}
		\end{subfigure}
		\hfill
		\begin{subfigure}{0.32\linewidth}
			\centering
			\includegraphics[width=\linewidth]{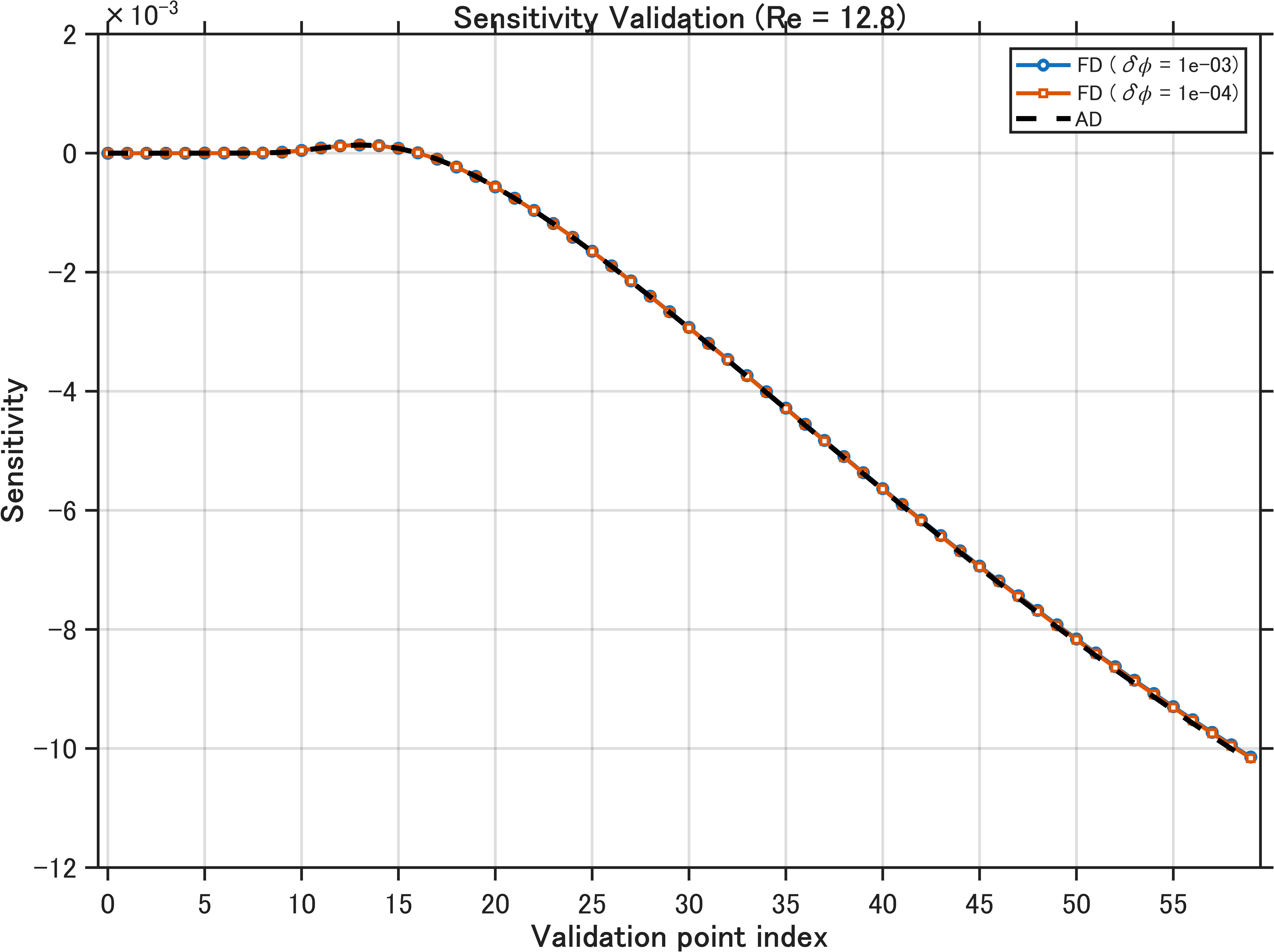}
			\caption{${\rm Re}=12.8$.}
			\label{sensitivityvalidationdifferentepsre12}
		\end{subfigure}
		\hfill
		\begin{subfigure}{0.32\linewidth}
			\centering
			\includegraphics[width=\linewidth]{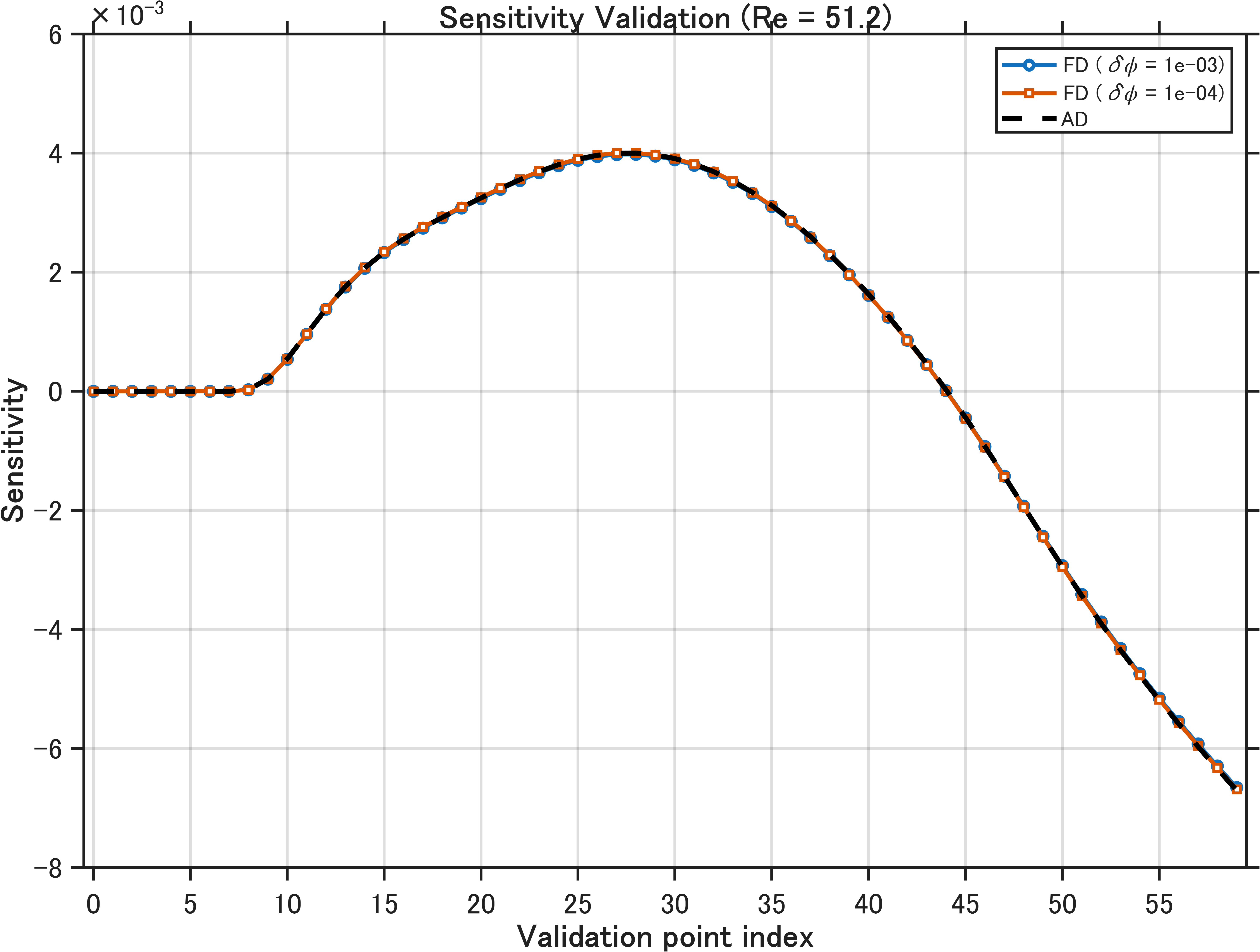}
			\caption{${\rm Re}=51.2$.}
			\label{sensitivityvalidationdifferentepsre51}
		\end{subfigure}
		
		\caption{Comparison between the finite-difference and adjoint sensitivities.}
		\label{sensitivityvalidation}
	\end{figure}

	The accuracy of the proposed adjoint sensitivity is first assessed by comparison with finite-difference sensitivities obtained from single-point design perturbations. For consistency, the finite-difference and adjoint sensitivities are evaluated using the same representative time window and sampling instants. Sixty perturbation points are uniformly selected downstream of the cylinder centre. At each point, two perturbation magnitudes, $10^{-3}$ and $10^{-4}$ are considered. Throughout the validation, a sampling interval of $\Delta W = 50$ and a representative time window length of 5000 time steps are employed. Fig.~\ref{sensitivityvalidation} compares the adjoint and finite-difference sensitivities for the three Reynolds numbers.
	
	As shown in Fig.~\ref{sensitivityvalidation}, the adjoint sensitivities closely reproduce the finite-difference sensitivity distributions for all three Reynolds numbers and both perturbation magnitudes. The corresponding relative errors remain below $1\%$, demonstrating that the proposed adjoint formulation accurately predicts the sensitivities in both steady and unsteady flow regimes.
	\subsubsection{Comparison between representative and prescribed time windows}
	To evaluate the optimization capability of the proposed framework, the vortex-shedding case at ${\rm Re}=51.2$ is considered and compared with the conventional prescribed-window strategy. The representative time windows identified by the proposed criterion contain 5000 time steps, and a sampling interval of $\Delta W=50$ is adopted. For the prescribed-time-window strategy, two time windows are considered: $t\in[20000,25000]$ and $t\in[5000,6000]$. The maximum number of optimization iterations is set to $N_{\max}=150$, and the optimization tolerance is prescribed as $\varepsilon_o=0.005$.
	
	\begin{figure}[htbp]
		\centering
		
		\begin{subfigure}{0.42\linewidth}
			\centering
			\includegraphics[width=\linewidth]{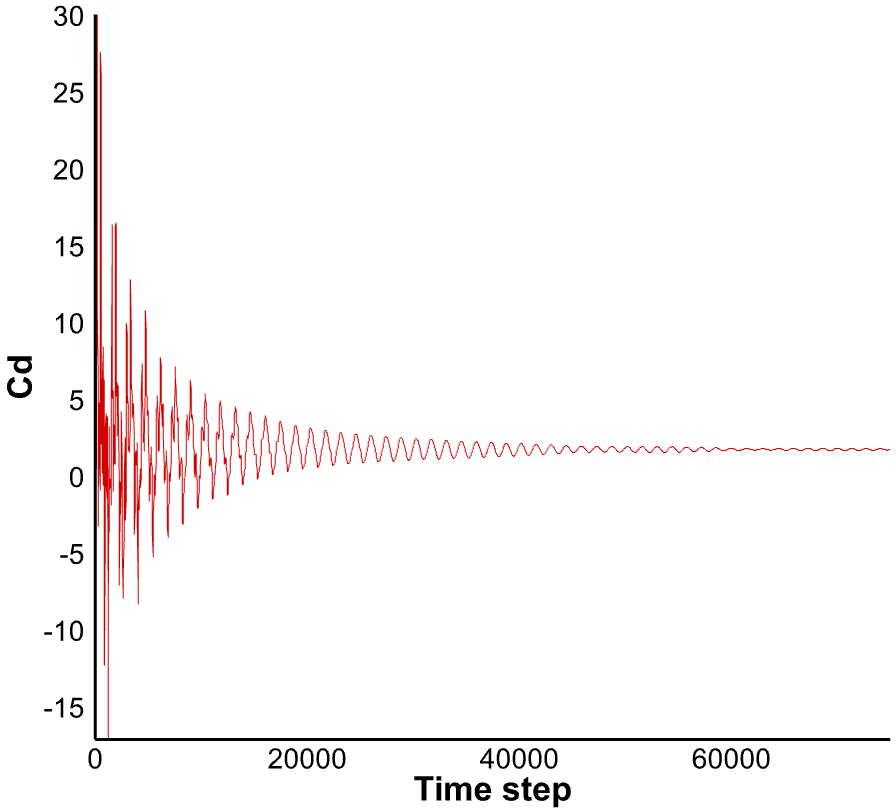}
			\caption{History of drag coefficient.}
			\label{cdhistory}
		\end{subfigure}
		\hfill
		\begin{subfigure}{0.48\linewidth}
			\centering
			\includegraphics[width=\linewidth]{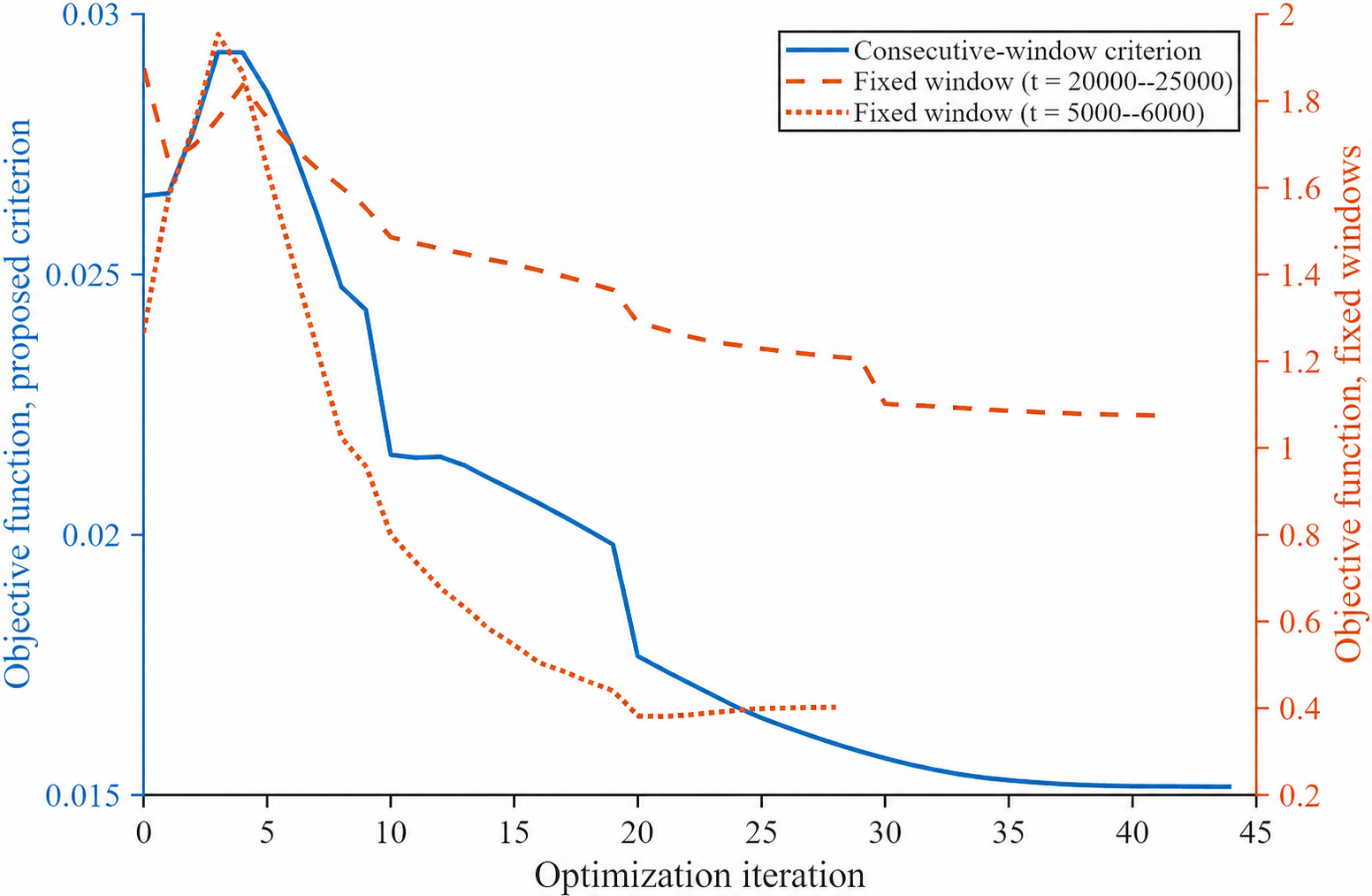}
			\caption{History of objective function.}
			\label{objectivefunctionhistorycriterion}
		\end{subfigure}
		
		\caption{Convergence histories.}
		\label{convergence_history}
	\end{figure}
	
\begin{figure}[htbp]
		\centering
		
		\begin{subfigure}{0.32\linewidth}
			\centering
			\includegraphics[width=\linewidth]{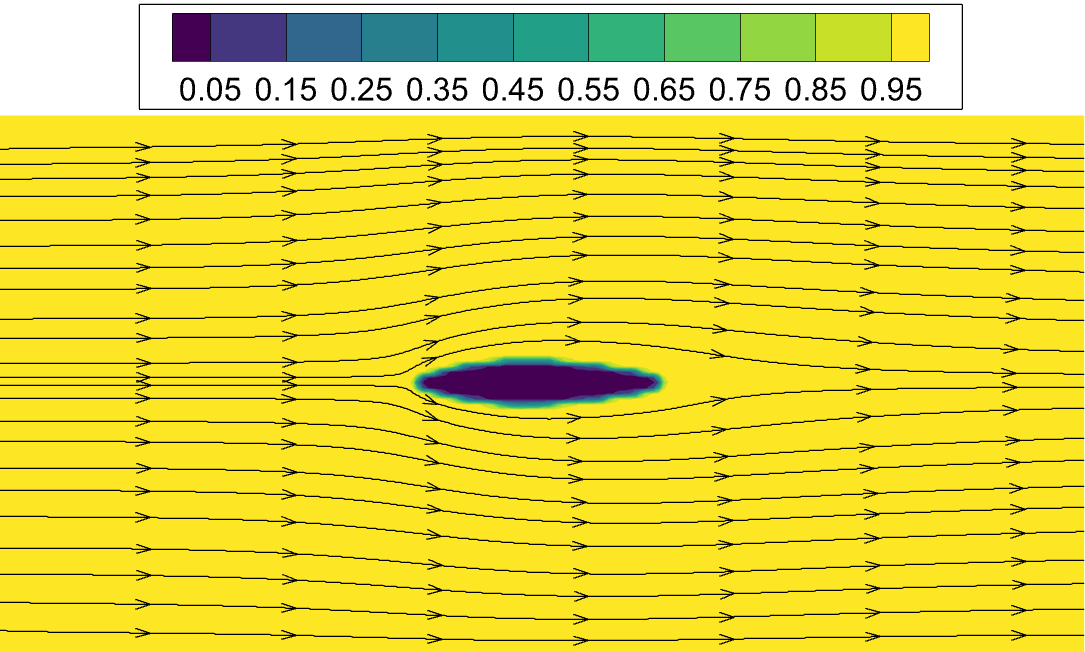}
			\caption{Proposed criterion.}
			\label{resultre51}
		\end{subfigure}
		\hfill
		\begin{subfigure}{0.32\linewidth}
			\centering
			\includegraphics[width=\linewidth]{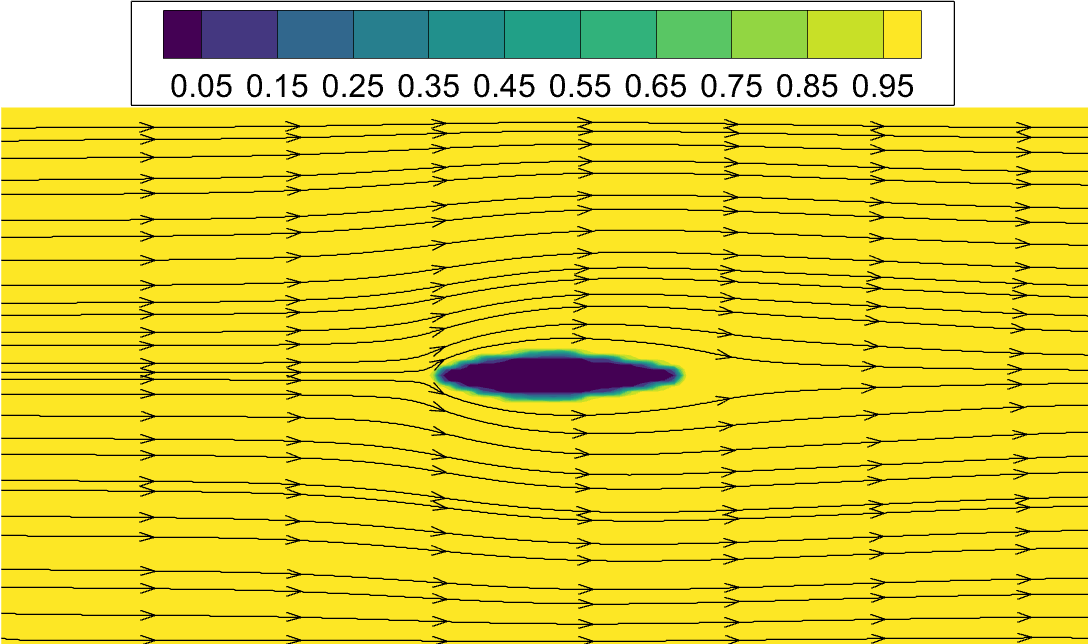}
			\caption{mild oscillations interval.}
			\label{resultmiddlewindow}
		\end{subfigure}
		\hfill
		\begin{subfigure}{0.32\linewidth}
			\centering
			\includegraphics[width=\linewidth]{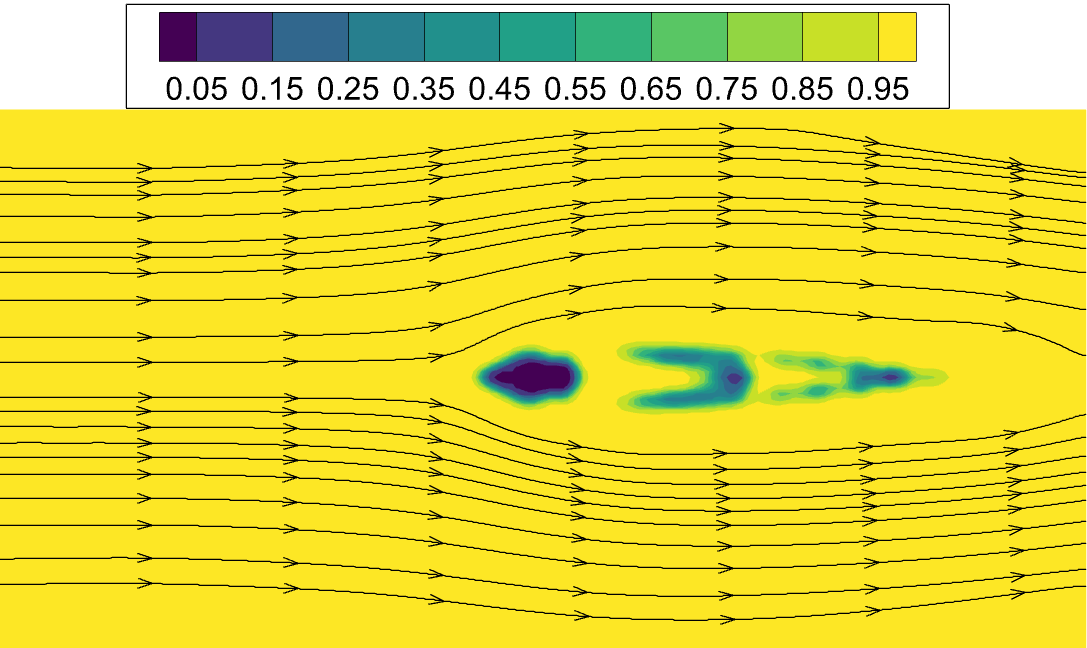}
			\caption{significant oscillations interval.}
			\label{resultearlywindow}
		\end{subfigure}
		
		\caption{Comparison of optimized results obtained using the proposed consecutive-window convergence criterion and the conventional fixed-window strategy.}
		\label{window_comparison}
	\end{figure}

	Figure~\ref{cdhistory} shows that the prescribed time windows correspond to different stages of flow development. During $t=5000$--$6000$, $C_d$ still exhibits pronounced early oscillations, whereas the oscillations become milder during $t=20000$--$25000$. Within the final time window identified by the consecutive-window convergence criterion, $C_d$ has approached an approximately constant level. The resulting optimized configurations are shown in Fig.~\ref{window_comparison}. The result obtained using the proposed criterion closely resembles that obtained using the later prescribed time window, and both recover the well-known streamlined, teardrop-shaped configuration. When reevaluated using the same representative-time-window criterion, their objective values differ by only approximately $0.1\%$, demonstrating the robustness of the optimization framework to the selection of a fully developed time window. In contrast, the early prescribed time window produces the physically unreasonable configuration shown in Fig.~\ref{resultearlywindow}.
	
	The corresponding objective-function histories are shown in Fig.~\ref{objectivefunctionhistorycriterion}. All three cases exhibit an overall decreasing trend. In particular, the objective associated with the early prescribed time window decreases despite producing an unreasonable configuration. Therefore, reduction of the objective function alone does not ensure that the prescribed time window represents the developed flow behavior.
	
	\subsubsection{Comparison with full-time adjoint optimization}
\begin{figure}[htbp]
		\centering
		
		\begin{subfigure}[b]{0.48\linewidth}
			\centering
			\includegraphics[width=\linewidth]{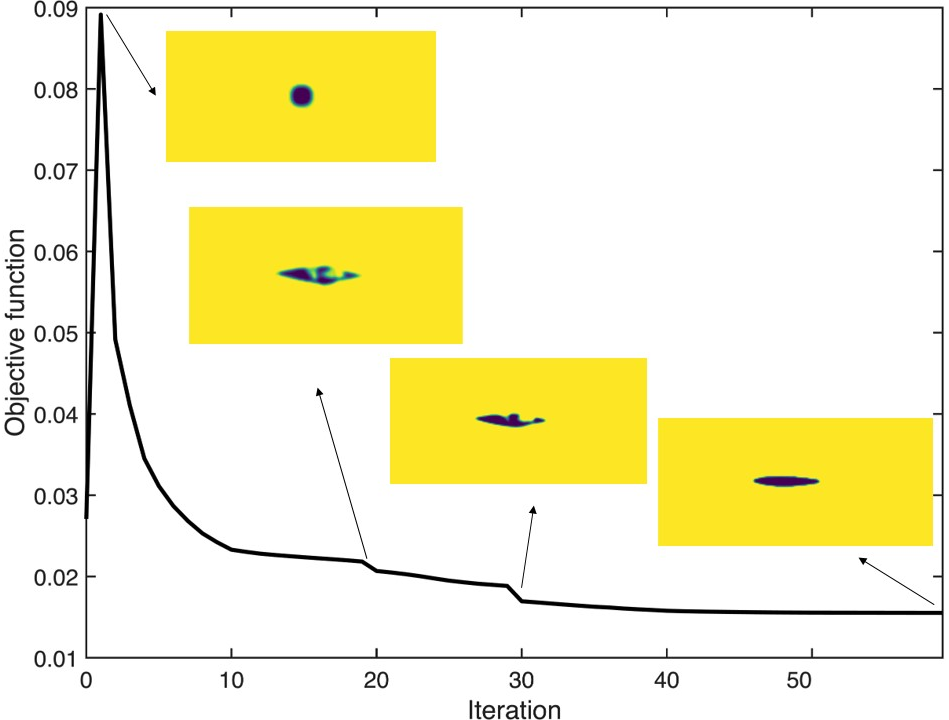}
			\caption{Full-time adjoint framework.}
			\label{fulltimeobj}
		\end{subfigure}
		\hfill
		\begin{subfigure}[b]{0.48\linewidth}
			\centering
			\includegraphics[width=\linewidth]{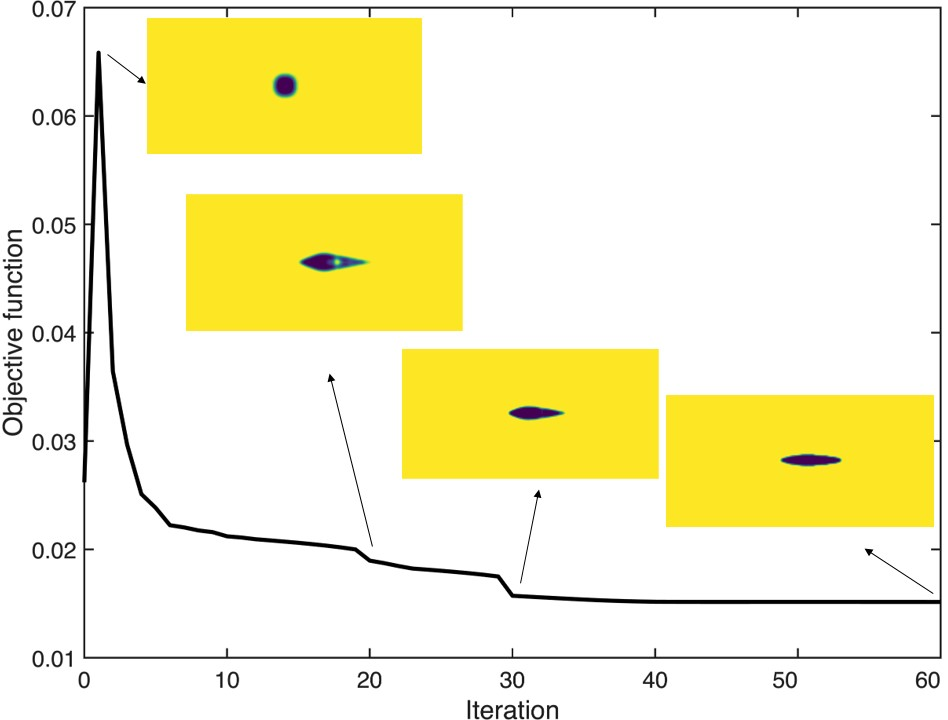}
			\caption{Proposed framework.}
			\label{60stepobj}
		\end{subfigure}
		
		\caption{Optimization histories obtained using full-time adjoint propagation and the proposed framework.}
		\label{fulltime_compar}
	\end{figure}

	To examine the influence of restricting the adjoint propagation to the representative time window, the cylinder-flow case at ${\rm Re}=51.2$ is considered. An additional optimization is performed using full-time adjoint propagation, in which the complete flow history is retained and the adjoint equations are traced backward to the initial time. For the proposed framework, the time-window length is set to 5000 time steps, with a sampling interval of $\Delta W=50$. Both approaches are performed for a fixed budget of $60$ optimization iterations. The corresponding objective-function histories and structural evolution are shown in Fig.~\ref{fulltime_compar}.
	
	Both final designs are reevaluated using the consecutive-window convergence criterion to ensure a consistent comparison. The resulting objective values differ by only approximately $0.003\%$, and the optimized configurations are nearly identical. These results demonstrate that restricting the adjoint propagation to the representative time window yields an optimized result equivalent to that obtained using full-time adjoint propagation.
	
	\begin{table}[htbp]
		\centering
		\begin{tabular}{lcc}
			\hline
			& Full-time framework & Proposed framework \\
			\hline
			Optimization iterations       & 60 & 60 \\
			Total computational time (min)  & 98 & 39.8 \\
			Average time per iteration (min)& 1.63 & 0.663 \\
			Storage requirement           & 558 GB & 217 MB \\
			Reevaluated objective value   & 0.01515151 & 0.01515106 \\
			\hline
		\end{tabular}
		\caption{Comparison of the computational costs and optimized results obtained using the proposed framework and full-time adjoint propagation.}
		\label{fulltime_comparison}
	\end{table}
	The average computational time per optimization iteration and the storage requirements are compared in Table~\ref{fulltime_comparison}. The proposed framework reduces the average runtime from $1.63$ to $0.663$ minutes per iteration, corresponding to a reduction of approximately $59.3\%$. It also reduces the required storage from $558$ GB to $217$ MB, a reduction of approximately $99.96\%$.
	
	\subsection{Wake flow recovery}
	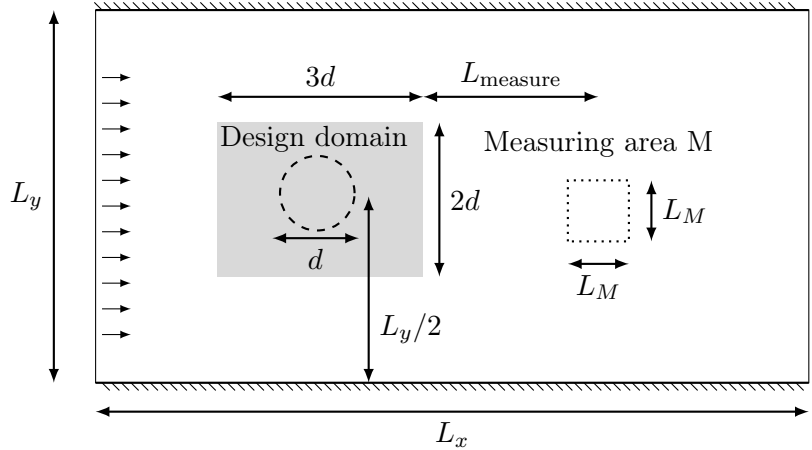
\begin{figure}[htbp]
		\centering
		\begin{tikzpicture}[
			scale=0.85,
			>=latex,
			font=\small
			]
			
			\def\Lx{12.0}
			\def\Ly{5.8}
			\def\xin{0.9}
			\def\yd{1.65}
			\def\h{2.4}
			\def\w{3.2}
			\def\xd{2.8}
			
			\draw (\xin,0) rectangle (\Lx,\Ly);
			
			\draw[thick] (\xin,\Ly) -- (\Lx,\Ly);
			\draw[thick] (\xin,0) -- (\Lx,0);
			
			\foreach \x in {1.0,1.15,...,11.8}{
				\draw[thin] (\x,\Ly) -- ++(-0.12,0.12);
				\draw[thin] (\x,0) -- ++(0.12,-0.12);
			}
			
			
			\foreach \y in {0.75,1.15,...,5.05}{
				\draw[->,thin] (1.0,\y) -- (1.45,\y);
			}
			
			\draw[<->,thick] (0.25,0.0) -- (0.25,\Ly);
			\node[left] at (0.25,\Ly/2) {$L_y$};
			
			\fill[gray!30] (\xd,\yd) rectangle ++(\w,\h);
			\node at (\xd+1.5,\yd+\h-0.25) {Design domain};
			
			\draw[dashed,thick] (\xd+1.55,\yd+1.30) circle (0.58);
			
			\draw[<->,thick] (\xd+0.85,\yd+0.60) -- ++(1.35,0);
			\node[below] at (\xd+1.52,\yd+0.60) {$d$};
			
			\draw[<->,thick] (\xd+\w+0.25,\yd) -- ++(0,\h);
			\node[right] at (\xd+\w+0.25,\yd+\h/2) {$2d$};
			
			\draw[<->,thick] (\xd+\w-0.85,0) -- (\xd+\w-0.85,\yd+1.25);
			\node[right] at (\xd+\w-0.85,0.85) {$L_y/2$};
			
			\draw[<->,thick] (\xd,\yd+\h+0.40) -- (\xd+\w,\yd+\h+0.40);
			\node[above] at (\xd+\w/2,\yd+\h+0.40) {$3d$};
			
			\def\xm{8.25}
			\def\ym{2.20}
			\def\Lm{0.95}
			
			\draw[dotted,thick] (\xm,\ym) rectangle ++(\Lm,\Lm);
			\node[above] at (\xm+\Lm/2,\ym+\Lm+0.20) {Measuring area M};
			
			\draw[<->,thick] (\xm,\ym-0.35) -- ++(\Lm,0);
			\node[below] at (\xm+\Lm/2,\ym-0.35) {$L_M$};
			
			\draw[<->,thick] (\xm+\Lm+0.35,\ym) -- ++(0,\Lm);
			\node[right] at (\xm+\Lm+0.35,\ym+\Lm/2) {$L_M$};
			
			\draw[<->,thick] (\xd+\w,\yd+\h+0.40) -- (\xm+\Lm/2,\yd+\h+0.40);
			\node[above] at ({(\xd+\w+\xm+\Lm/2)/2},\yd+\h+0.40) {$L_{\rm measure}$};

			\draw[<->,thick] (\xin,-0.45) -- (\Lx,-0.45);
			\node[below] at ({(\xin+\Lx)/2},-0.45) {$L_x$};
			
		\end{tikzpicture}
		
		\caption{Illustration of the computational domain for wake flow recovery problem.}
		\label{wake_flow_recovery}
	\end{figure}
	To validate the capability of the proposed framework to capture unsteady flow characteristics, the wake flow recovery problem is considered as a benchmark. In this example, the objective is to recover the wake generated by an obstacle through the optimization of its structural configuration. 
	\begin{table}[htbp]
		\centering
		\begin{tabular}{lc}
			\hline
			Parameter & Value \\
			\hline
			$L_x$ & 200 \\
			$L_y$ & 50 \\
			$d$ & 10, 12 \\
			$L_M$ & 6 \\
			$L_{\rm measure}$ & 33, 39 \\
			$|W_i|$ & 6000 \\
			$\Delta W$ & 20 \\
			$\varepsilon$ & $5\times10^{-4}$ \\
			$\varepsilon_o$ & $5\times10^{-3}$ \\
			$N_{\max}$ & 200 \\
			\hline
		\end{tabular}
		\caption{Computational parameters for the wake flow recovery problem.}
		\label{wfr}
	\end{table}
	As illustrated in Fig.~\ref{wake_flow_recovery}, a circular cylinder is first placed in the design domain $D$, and the corresponding time-averaged flow field $\mathbf{\bar{u}}_{\rm ref}$ is recorded after the consecutive-window convergence criterion is satisfied. Subsequently, the initial material distribution in the design domain is set to a constant value within the prescribed volume constraint, where the allowable volume fraction ranges from $0.1$ to $0.9$. During the optimization process, the objective is to minimize the discrepancy between the time-averaged velocity field within the measurement domain $M$ and the reference flow field, which is mathematically expressed as
	\begin{equation}
		\begin{aligned}
			&\min_{\phi} J=\sum_{i\in M}|\mathbf{\bar{u}}(\mathbf{x}_i)-\mathbf{\bar{u}}_{\rm ref}(\mathbf{x}_i)|^2,\\
			&\text{s.t.}\quad 0.1\leq\frac{1}{n_D}\sum_{i\in D}(1-\phi_i) \leq0.9,
		\end{aligned}
	\end{equation}
	where 
	\begin{equation}
		\mathbf{\bar{u}}(\mathbf{x}_i)=\sum_{t_n \in W_{\rm eff}}\frac{\Delta W}{|W_{\rm eff}|}\mathbf{u}(\mathbf{x}_i,t_n).
	\end{equation}
	Here, $\mathbf{\bar{u}}(\mathbf{x}_i)$ denotes the time-averaged velocity field during optimization process.
	The corresponding computational parameters are summarized in Table~\ref{wfr}. The values of $d$ and $L_{\rm measure}$ are specified separately for the ${\rm Re}=60$ and ${\rm Re}=72$ cases.
	\begin{figure}[htbp]
		\centering
		
		\begin{subfigure}{0.48\linewidth}
			\centering
			\includegraphics[width=\linewidth]{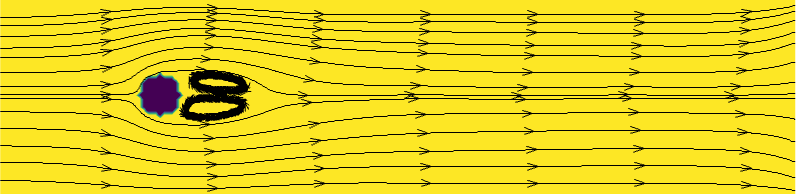}
			\caption{Reference design.}
			\label{wakeflowinitialre60phi}
		\end{subfigure}
		\hfill
		\begin{subfigure}{0.48\linewidth}
			\centering
			\includegraphics[width=\linewidth]{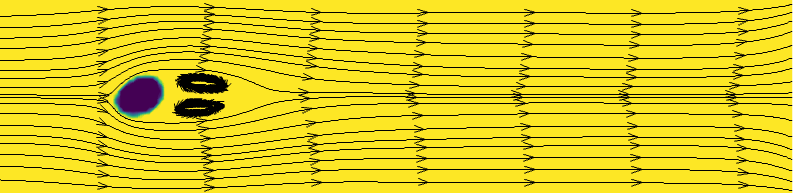}
			\caption{Optimized design.}
			\label{wakeflowresultre60phi}
		\end{subfigure}
		
		\vspace{0.5em}
		
		\begin{subfigure}{0.48\linewidth}
			\centering
			\includegraphics[width=\linewidth]{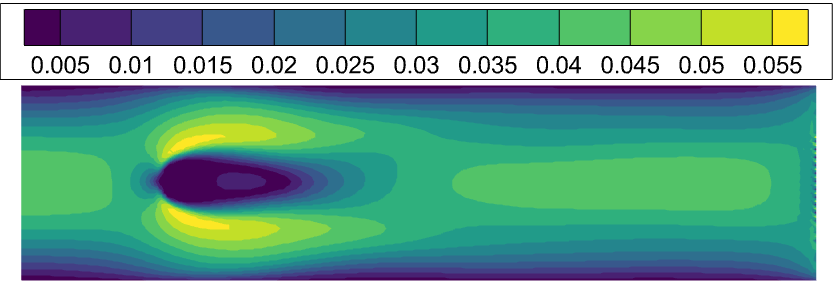}
			\caption{Reference velocity field.}
			\label{wakeflowinitialre60velocity}
		\end{subfigure}
		\hfill
		\begin{subfigure}{0.48\linewidth}
			\centering
			\includegraphics[width=\linewidth]{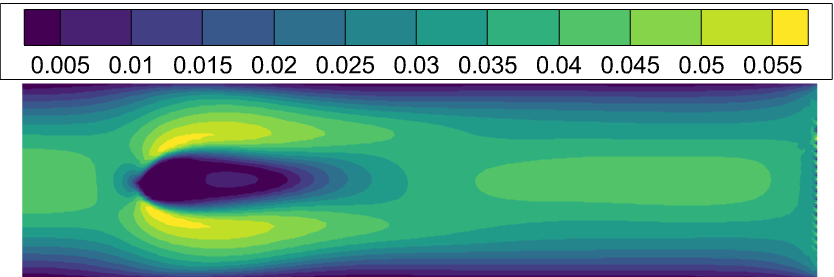}
			\caption{Optimized velocity field.}
			\label{wakeflowresultre60velocity}
		\end{subfigure}
		
		\caption{Reference and optimized configurations for the wake flow recovery problem at ${\rm Re}=60$.}
		\label{wakeflowre60}
	\end{figure}
	
\begin{figure}[htbp]
		\centering
		
		\begin{subfigure}{0.48\linewidth}
			\centering
			\includegraphics[width=\linewidth]{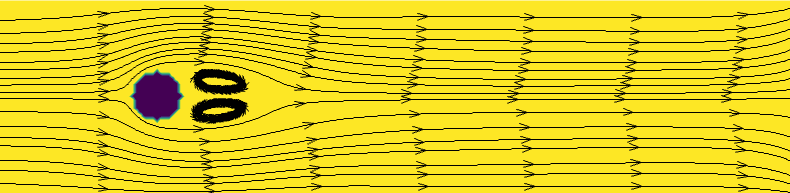}
			\caption{Reference design.}
			\label{wakeflowinitialre72phi}
		\end{subfigure}
		\hfill
		\begin{subfigure}{0.48\linewidth}
			\centering
			\includegraphics[width=\linewidth]{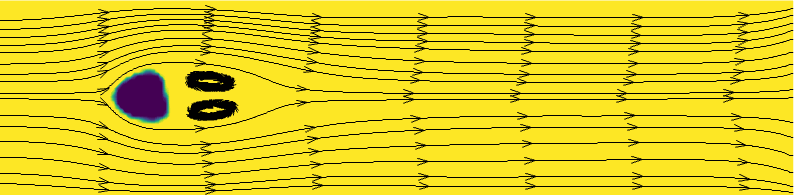}
			\caption{Optimized design.}
			\label{wakeflowresultre72phi}
		\end{subfigure}
		
		\vspace{0.5em}
		
		\begin{subfigure}{0.48\linewidth}
			\centering
			\includegraphics[width=\linewidth]{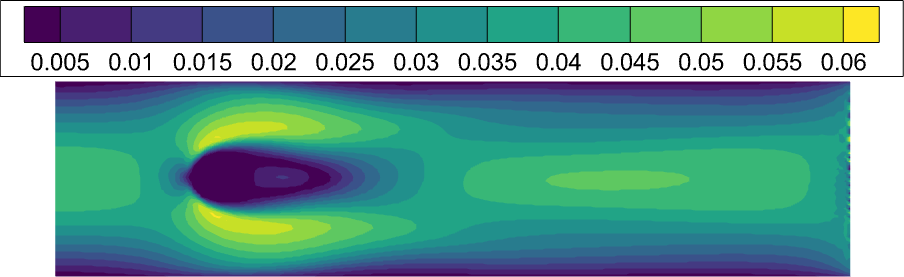}
			\caption{Reference velocity field.}
			\label{wakeflowinitialre72velocity}
		\end{subfigure}
		\hfill
		\begin{subfigure}{0.48\linewidth}
			\centering
			\includegraphics[width=\linewidth]{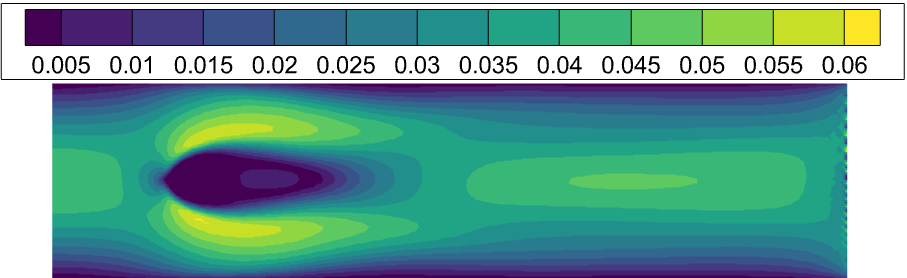}
			\caption{Optimized velocity field.}
			\label{wakeflowresultre72velocity}
		\end{subfigure}
		
		\caption{Reference and optimized configurations for the wake flow recovery problem at ${\rm Re}=72$.}
		\label{wakeflowre72}
	\end{figure}

	For the cases with ${\rm Re}=60$ and ${\rm Re}=72$, the reference time-averaged flow fields and the corresponding optimized configurations are presented in Fig.~\ref{wakeflowre60} and~\ref{wakeflowre72}. It can be observed that, although vortex shedding occurs in both cases, the time-averaged flow field exhibits a pair of symmetric recirculation regions. For both Reynolds numbers, the optimized structures are successfully recovered to an approximately circular-cylinder shape. 
	
	For ${\rm Re}=60$, the objective function decreases from $3.064\times10^{-5}$ to $2.02\times10^{-7}$, while for ${\rm Re}=72$, it decreases from $3.94\times10^{-5}$ to $7.77\times10^{-8}$. In both cases, the objective function is reduced by more than $99\%$. These results indicate that, for unsteady flows, the time-averaged state still preserves the essential characteristics of the unsteady flow, while the proposed criterion correctly determines whether the flow evolution has reached a fully developed stage.
	
	\subsection{U-bend}
	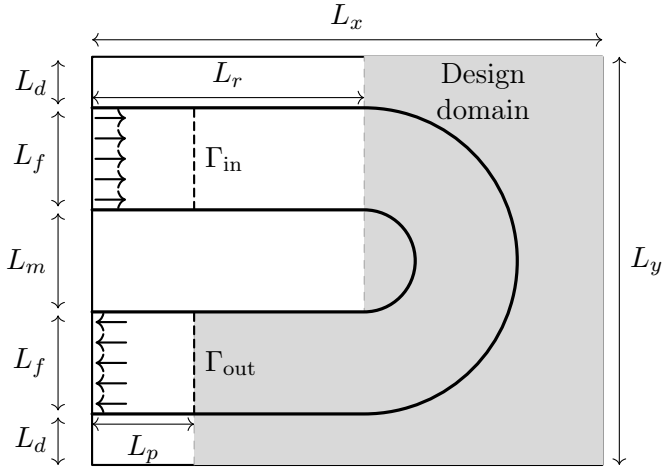
\begin{figure}[htbp]
		\centering
		\begin{tikzpicture}[
			x=4.5cm,
			y=4.5cm,
			line cap=round,
			line join=round,
			font=\small
			]
			\def\Lx{1.5}
			\def\Ly{1.2}
			
			\def\xc{0.8}
			\def\yc{0.6}
			
			\def\Rout{0.45}
			\def\Rin{0.15}
			
			\def\ylowmin{0.15}
			\def\ylowmax{0.45}
			\def\yupmin{0.75}
			\def\yupmax{1.05}
			
			\def\xdesign{0.3}
			
			\draw[thick] (0,0) rectangle (\Lx,\Ly);
			
			\fill[gray!30]
			(\xc,0) rectangle (\Lx,\Ly);
			
			\fill[gray!30]
			(\xdesign,0) rectangle (\xc,\ylowmax);
			
			\draw[gray!70, dashed]
			(\xdesign,0)
			-- (\xdesign,\ylowmax)
			-- (\xc,\ylowmax)
			-- (\xc,\Ly);
			
			\node[align=center] at (1.15,1.10)
			{Design\\domain};
			\draw[very thick]
			(0,\yupmax)
			-- (\xc,\yupmax)
			arc[
			start angle=90,
			end angle=-90,
			radius=\Rout
			]
			-- (0,\ylowmin);
			
			\draw[very thick]
			(0,\yupmin)
			-- (\xc,\yupmin)
			arc[
			start angle=90,
			end angle=-90,
			radius=\Rin
			]
			-- (0,\ylowmax);
			
			\foreach \y in {0.78,0.84,0.90,0.96,1.02}{
				\draw[->, thick] (0.01,\y) -- (0.10,\y);
			}
			\foreach \y in {0.18,0.24,0.30,0.36,0.42}{\draw[->, thick] (0.10,\y) -- (0.01,\y);}
			\draw[<->] (-0.1,0) -- (-0.1,0.14);
			\node[left] at (-0.1,0.075) {$L_d$};
			
			\draw[<->] (-0.1,0.16) -- (-0.1,0.44);
			\node[left] at (-0.1,0.30) {$L_f$};
			
			\draw[<->] (-0.1,0.46) -- (-0.1,0.74);
			\node[left] at (-0.1,0.60) {$L_m$};
			
			\draw[<->] (-0.1,0.76) -- (-0.1,1.04);
			\node[left] at (-0.1,0.90) {$L_f$};
			
			\draw[<->] (-0.1,1.06) -- (-0.1,1.19);
			\node[left] at (-0.1,1.125) {$L_d$};
			
			\draw[<->] (1.55,0) -- (1.55,\Ly);
			\node[right] at (1.55,0.60) {$L_y$};
			\draw[<->] (0,0.12) -- (\xdesign,0.12);
			\node[below] at (0.15,0.12) {$L_p$};
			
			\draw[<->] (0,1.08) -- (\xc,1.08);
			\node[above] at (0.4,1.08) {$L_r$};
			
			\draw[<->] (0,1.25) -- (\Lx,1.25);
			\node[above] at (0.75,1.25) {$L_x$};
			\def\xpressure{0.3}
			\draw[densely dashed, thick]
			(\xpressure,\yupmin) -- (\xpressure,\yupmax);
			\draw[densely dashed, thick]
			(\xpressure,\ylowmin) -- (\xpressure,\ylowmax);
			\node[right] at (\xpressure,0.9)
			{$\Gamma_{\mathrm{in}}$};
			\node[right] at (\xpressure,0.3)
			{$\Gamma_{\mathrm{out}}$};
		\end{tikzpicture}
		
		\caption{Schematic of the design domain for U-bend flow.}
		\label{ubend_geometry}
	\end{figure}
	
	Next, the proposed framework is applied to a higher-Reynolds-number U-bend flow. In this configuration, a uniform inlet velocity enters the design domain and passes through a $180^\circ$ sharp bend. As the flow turns into the downstream channel, flow separation and vortex formation readily occur, resulting in pronounced unsteady flow behavior. 
	
	The computational configuration is illustrated in Fig.~\ref{ubend_geometry}, where $L_f$ denotes the channel width, $L_m$ represents the diameter of the turning section connecting the upper and lower channels, and $L_d$ is the distance between the channel wall and the computational boundary. The quantities $L_p$ and $L_r$ denote the reserved outlet and inlet channel lengths, respectively. The gray region indicates the design domain. 
	\begin{table}[htbp]
		\centering
		\begin{tabular}{lcc}
			\toprule
			Parameter & Re = 400 & Re = 750 \\
			\midrule
			$L_x$               & 240   & 300 \\
			$L_y$               & 192   & 240 \\
			$L_d$               & 24    & 30 \\
			$L_f$               & 48    & 60 \\
			$L_m$               & 48    & 60 \\
			$L_p$               & 48    & 60 \\
			$L_r$               & 144   & 180 \\
			$|W_i|$             & 20000 & 20000 \\
			$\Delta W$        & 20    & 20 \\
			$\varepsilon$              & 0.05   & 0.05 \\
			$\varepsilon_o$       & $10^{-5}$ & $10^{-5}$ \\
			$N_{\max}$          & 90   & 90 \\
			\bottomrule
		\end{tabular}
		\caption{Computational parameters for the U-bend optimization cases.}
		\label{ubend_parameters}
	\end{table}
	
	In this configuration, the objective is to minimize the pressure drop across the channel. The inlet and outlet pressure measurement sections are denoted by $\Gamma_{\mathrm{in}}$ and $\Gamma_{\mathrm{out}}$, respectively. The objective function is mathematically expressed as
	\begin{equation}
		\begin{aligned}
			&\min_{\phi} J=\frac{\Delta W}{|W_{\rm eff}|}\sum_{t_n\in W_{\rm eff}}\left(\bar{p}_{\rm in}(t_n)-\bar{p}_{\rm out}(t_n)\right),\\
			&\text{s.t.}\quad 0.52\leq\frac{1}{n_D}\sum_{i\in D}(1-\phi_i) \leq0.65,
		\end{aligned}
	\end{equation}
	where $\bar{p}_{\rm in}$ and $\bar{p}_{\rm out}$ denote the mean pressures over the inlet and outlet measurement sections, respectively.
	
	Two representative cases with Reynolds numbers of $\rm Re = 400$ and $\rm Re = 750$ are considered. The corresponding computational parameters are summarized in Table~\ref{ubend_parameters}.
	
\begin{figure}[htbp]
		\centering
		
		\begin{subfigure}[b]{0.4\linewidth}
			\centering
			\includegraphics[width=\linewidth]{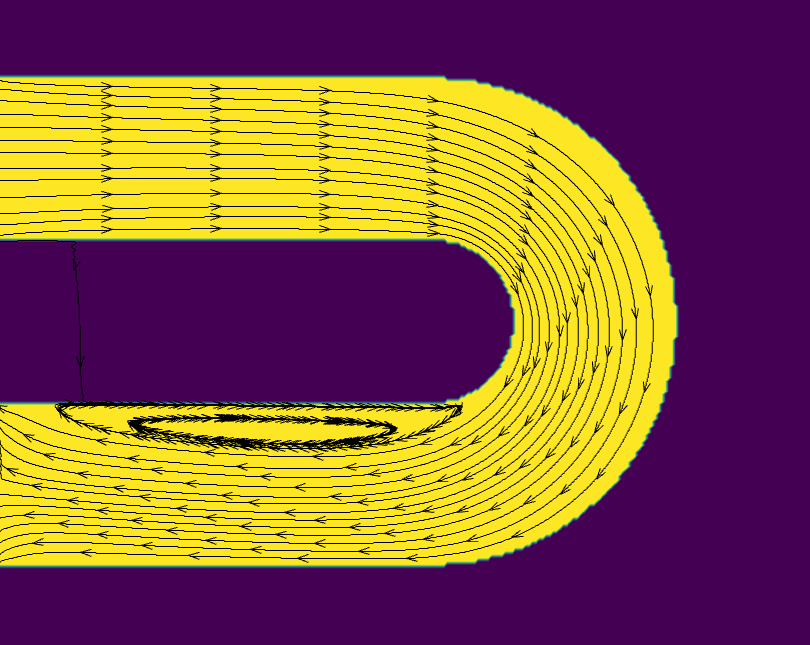}
			\caption{Re = 400}
		\end{subfigure}
		\hfill
		\begin{subfigure}[b]{0.4\linewidth}
			\centering
			\includegraphics[width=\linewidth]{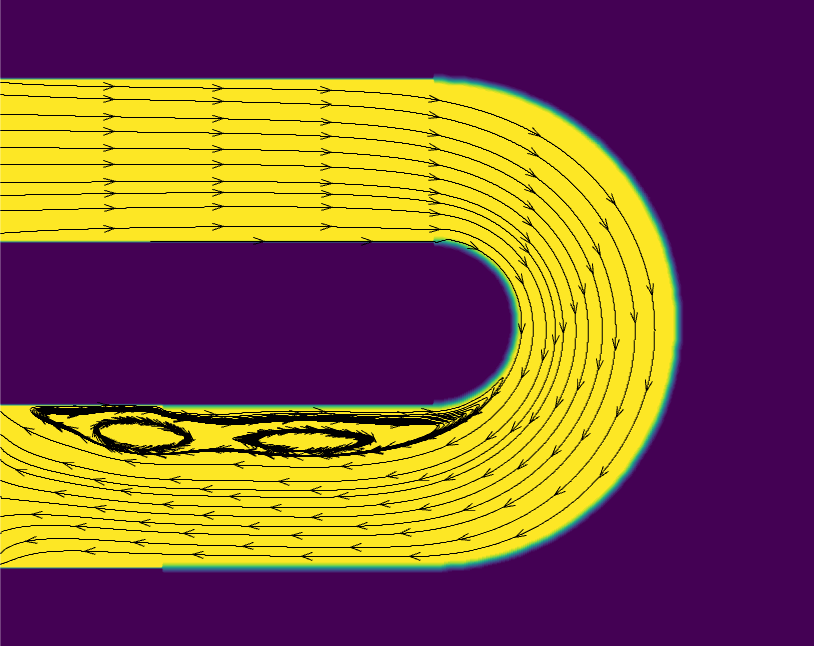}
			\caption{Re = 750}
		\end{subfigure}
		
		\caption{Initial U-bend configurations and instantaneous flow fields corresponding to Fig.~4 of Zhang and Pothérat~\cite{Zhang2013}.}
		\label{ubend_initial}
	\end{figure}

	According to the classification reported by Zhang and Pothérat~\cite{Zhang2013}, these two cases represent distinct unsteady flow regimes. Specifically, $\rm Re = 400$ corresponds to the onset of periodic vortex shedding in the wake region, whereas at $\rm Re = 750$, the downstream recirculation region has already broken down, giving rise to a more complex unsteady flow pattern with fragmented vortical structures. The corresponding initial configurations and instantaneous flow fields are shown in Fig.~\ref{ubend_initial}.
	
\begin{figure}[htbp]
		\centering
		
		\begin{subfigure}[b]{0.48\linewidth}
			\centering
			\includegraphics[width=\linewidth]{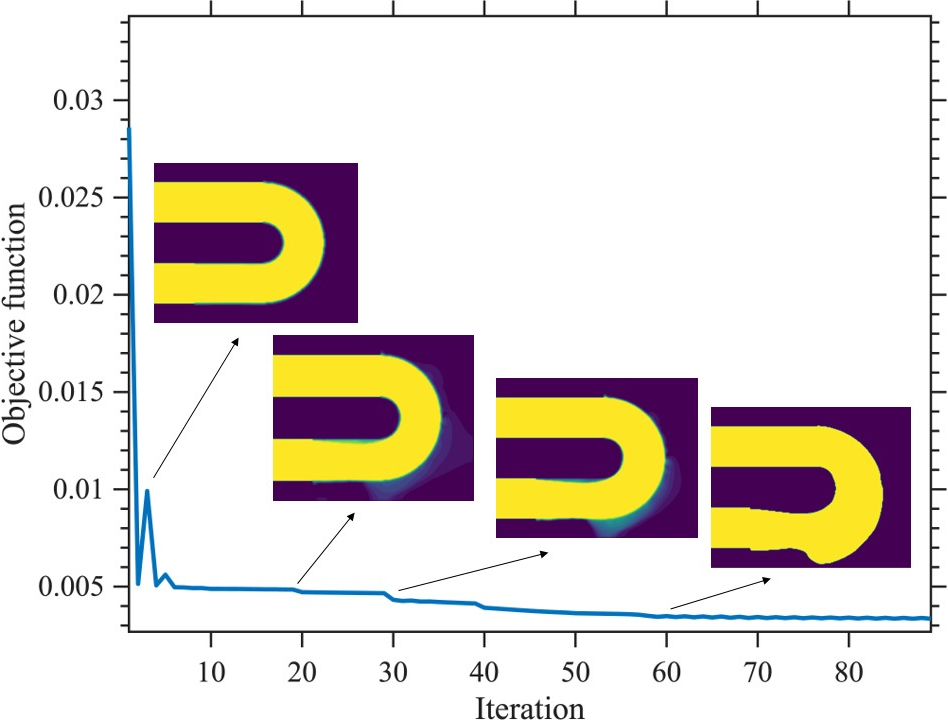}
			\caption{Re = 400}
		\end{subfigure}
		\hfill
		\begin{subfigure}[b]{0.48\linewidth}
			\centering
			\includegraphics[width=\linewidth]{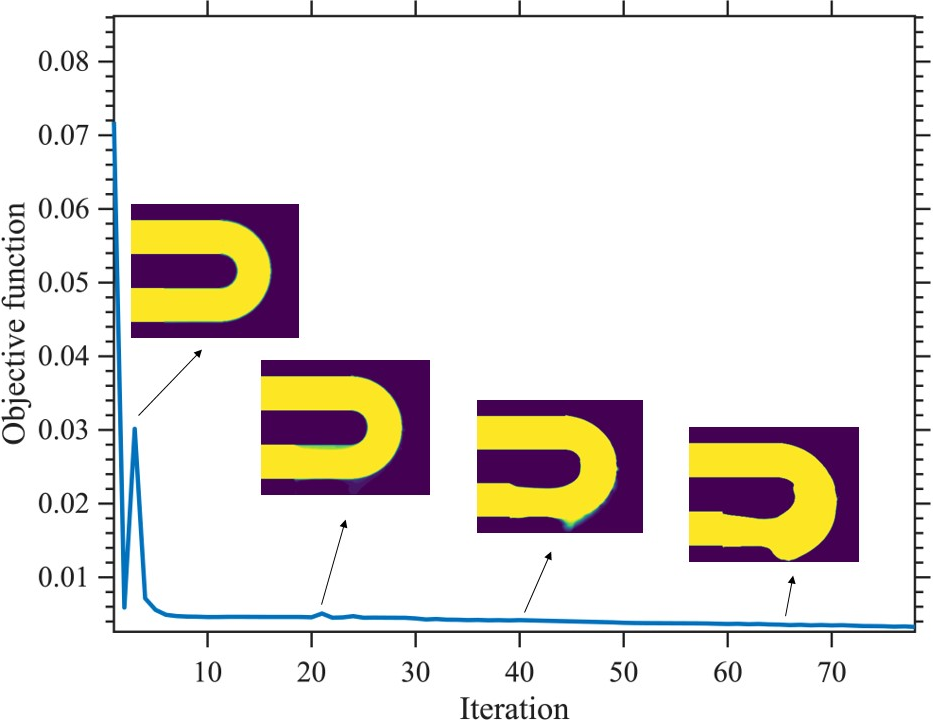}
			\caption{Re = 750}
		\end{subfigure}
		
		\caption{Evolution histories of the objective function for the U-bend optimization at $\mathrm{Re}=400$ and $\mathrm{Re}=750$.}
		\label{ubend_objective}
	\end{figure}

	The corresponding objective-function histories are shown in Fig.~\ref{ubend_objective}. As the optimization proceeds, the channel geometry initially evolves to suppress the large downstream recirculation region while enlarging the downstream buffer region, thereby mitigating the unsteady vortical structures induced by flow separation and backflow. Subsequently, the turning section expands, providing additional space to accommodate the change in flow direction and ultimately giving rise to a localized recirculation region near the downstream inner wall. During the later stages of optimization, particularly for the more unstable case at $\rm Re=750$, a vertically oriented recirculation region develops near the inner wall. This recirculation structure facilitates the discharge of the main flow and consequently reduces the pressure drop. Overall, the objective function decreases from $0.0286037$ to $0.00338550$ at $\rm Re=400$, corresponding to a reduction of $88.16\%$, and from $0.0718148$ to $0.00433962$ at $\rm Re=750$, corresponding to a reduction of $93.96\%$.
	
	\begin{figure}[htbp]
		\centering
		
		\begin{subfigure}[b]{0.365\linewidth}
			\centering
			\includegraphics[width=\linewidth]{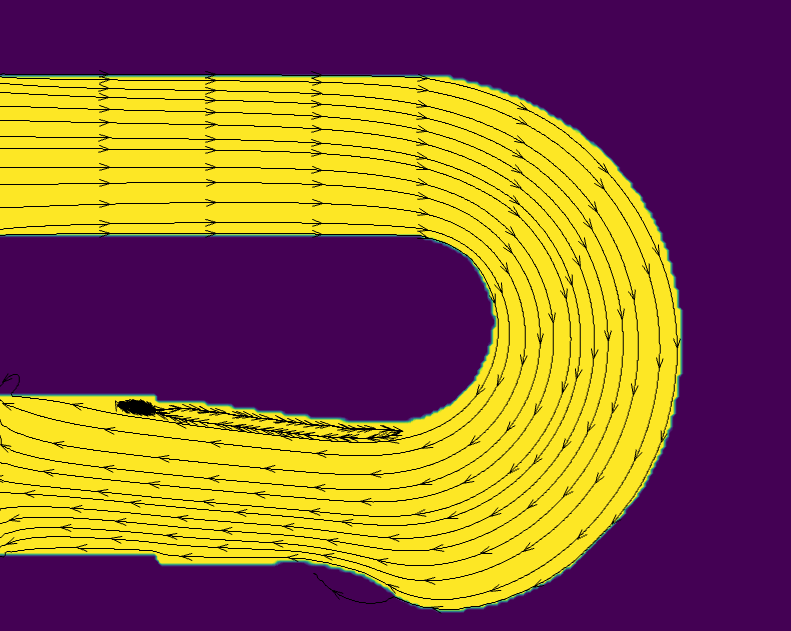}
			\caption{Re = 400}
		\end{subfigure}
		\hfill
		\begin{subfigure}[b]{0.365\linewidth}
			\centering
			\includegraphics[width=\linewidth]{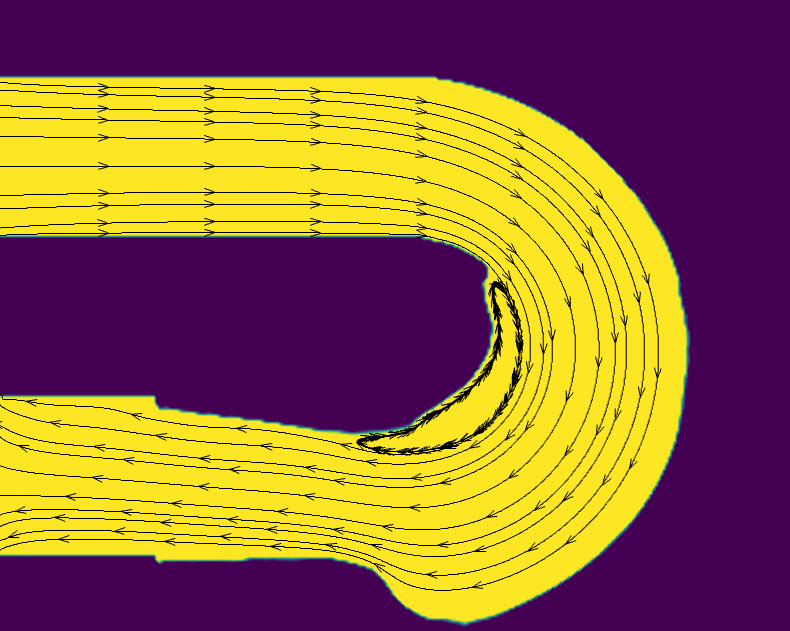}
			\caption{Re = 750}
		\end{subfigure}
		
		\caption{Optimized configurations and the corresponding instantaneous flow fields for the U-bend cases at $\mathrm{Re}=400$ and $\mathrm{Re}=750$.}
		\label{ubend_result}
	\end{figure}
	
	The optimized configurations and their corresponding flow fields are presented in Fig.~\ref{ubend_result}. Compared with the case at ${\rm Re}=400$, the optimized geometry at the higher Reynolds number, ${\rm Re}=750$, exhibits a larger buffer region downstream of the bend, providing additional space to accommodate the stronger flow turning and associated unsteadiness. Meanwhile, a localized, vertically elongated recirculation structure develops near the inner wall of the turning section. This structure appears to assist in redirecting the main flow through the bend. These results suggest that the proposed framework does not reduce the pressure drop by indiscriminately suppressing all vortical motions. Instead, it reorganizes the flow field by suppressing the large downstream recirculation structures responsible for substantial pressure losses while preserving or promoting localized structures that are potentially beneficial to flow transport. The U-bend cases therefore demonstrate the capability of the proposed framework to optimize unsteady flows at relatively high Reynolds numbers.
	
	\FloatBarrier
	\section{Conclusions}\label{Conclusion}
	In this work, an adjoint topology optimization framework with representative time windows has been developed for unsteady incompressible flows. By applying the consecutive-window convergence criterion during each forward solve, the framework automatically identifies a representative time window for every evolving structure, without prescribing its location in advance. The objective function and backward discrete adjoint analysis are then consistently restricted to the identified time window.
	The proposed framework was validated through a series of benchmark and topology optimization problems. The backward-facing-step problem first confirmed the accuracy of the forward LBM-LES solver, while the cylinder-flow cases demonstrated that the discrete adjoint sensitivities agreed well with the corresponding finite-difference results across attached, steady separated, and periodically unsteady flow regimes. The wake-flow recovery problems further verified the capability of the representative-window strategy to capture the essential time-averaged characteristics of unsteady flows and recover the prescribed reference configurations. Finally, the U-bend optimization cases demonstrated that the framework can accommodate more complex unsteady flow patterns and reduce pressure loss by reorganizing, rather than simply suppressing, the vortical structures. Taken together, these results demonstrate the effectiveness of the proposed framework in adaptively identifying representative time windows for evolving designs and producing physically meaningful optimized configurations under the flow regimes considered.
	
	Future work will focus on extending the present framework to three-dimensional turbulent flows, where more complex spatial and temporal dynamics must be considered. The framework will also be extended to coupled thermal-flow topology optimization problems, enabling the simultaneous treatment of unsteady fluid transport and heat-transfer performance in more complex engineering applications.
	
	\FloatBarrier
	\section*{Code and data availability}
	
	The code and data supporting the findings of this study are available from the authors upon reasonable request.
	
	\section*{Funding}
	
	This work was supported in part by the National Key Basic
	Research Program [grant number 2022YFA1004402]; the National
	Natural Science Foundation of China [grant number 12471377];
	and the Science and Technology Commission of Shanghai
	Municipality [grant number 22DZ2229014].
	\bibliographystyle{elsarticle-num}
	\bibliography{reference}

@article{Junk2005,
	author  = {Junk, M. and Klar, A. and Luo, L.-S.},
	title   = {Asymptotic Analysis of the Lattice {Boltzmann} Equation},
	journal = {J. Comput. Phys.},
	year    = {2005},
	volume  = {210},
	number  = {2},
	pages   = {676--704},
	doi     = {10.1016/j.jcp.2005.05.003}
}

@article{Smagorinsky1963,
	author  = {Smagorinsky, J.},
	title   = {General Circulation Experiments with the Primitive Equations: I. The Basic Experiment},
	journal = {Mon. Weather Rev.},
	year    = {1963},
	volume  = {91},
	number  = {3},
	pages   = {99--164}
}

@article{Lallemand2000,
	author  = {Lallemand, P. and Luo, L.-S.},
	title   = {Theory of the Lattice {Boltzmann} Method: Dispersion, Dissipation, Isotropy, Galilean Invariance, and Stability},
	journal = {Phys. Rev. E},
	year    = {2000},
	volume  = {61},
	number  = {6},
	pages   = {6546--6562},
	doi     = {10.1103/PhysRevE.61.6546}
}

@incollection{Hou1996,
	author    = {S. Hou and J. Sterling and S. Chen and G. D. Doolen},
	title     = {A Lattice {Boltzmann} Subgrid Model for High Reynolds Number Flows},
	booktitle = {Pattern Formation and Lattice Gas Automata},
	editor    = {A. T. Lawniczak and R. Kapral},
	series    = {Fields Inst. Commun.},
	volume    = {6},
	pages     = {151--166},
	year      = {1996},
	publisher = {Am. Math. Soc.},
	address   = {Providence, RI}
}

@article{Latt2006,
	author  = {Latt, J. and Chopard, B.},
	title   = {Lattice {Boltzmann} Method with Regularized Pre-Collision Distribution Functions},
	journal = {Math. Comput. Simul.},
	year    = {2006},
	volume  = {72},
	number  = {2--6},
	pages   = {165--168},
	doi     = {10.1016/j.matcom.2006.05.017}
}

@article{Zou1997,
	author  = {Zou, Q. and He, X.},
	title   = {On Pressure and Velocity Boundary Conditions for the Lattice {Boltzmann BGK} Model},
	journal = {Phys. Fluids},
	year    = {1997},
	volume  = {9},
	number  = {6},
	pages   = {1591--1598},
	doi     = {10.1063/1.869307}
}

@article{Borrvall2003,
	author  = {Borrvall, T. and Petersson, J.},
	title   = {Topology Optimization of Fluids in {Stokes} Flow},
	journal = {Int. J. Numer. Methods Fluids},
	year    = {2003},
	volume  = {41},
	number  = {1},
	pages   = {77--107},
	doi     = {10.1002/fld.426}
}

@article{Norgaard2016,
	author  = {N{\o}rgaard, S. A. and Sigmund, O. and Lazarov, B. S.},
	title   = {Topology Optimization of Unsteady Flow Problems Using the Lattice {Boltzmann} Method},
	journal = {J. Comput. Phys.},
	year    = {2016},
	volume  = {307},
	pages   = {291--307},
	doi     = {10.1016/j.jcp.2015.12.023}
}

@article{Zhu2013,
	author  = {J. Zhu and J. Ma},
	title   = {An Improved Gray Lattice {Boltzmann} Model for Simulating Fluid Flow in Multi-Scale Porous Media},
	journal = {Adv. Water Resour.},
	year    = {2013},
	volume  = {56},
	pages   = {61--76},
	issn    = {0309-1708},
	doi     = {10.1016/j.advwatres.2013.03.001}
}

@article{Ramalingom2025,
	author  = {Delphine Ramalingom and Alain Bastide and
	Pierre-Henri Cocquet and Micha{\"e}l Rakotob{\'e} and
	David Marti},
	title   = {Time-averaged algorithm for solving the topology
	optimization problem for unsteady laminar, turbulent
	and anisothermal flows},
	journal = {Journal of Computational Physics},
	volume  = {538},
	pages   = {114175},
	year    = {2025},
	doi     = {10.1016/j.jcp.2025.114175}
}

@article{Williamson1996,
	author  = {Williamson, C. H. K.},
	title   = {Vortex Dynamics in the Cylinder Wake},
	journal = {Annu. Rev. Fluid Mech.},
	year    = {1996},
	volume  = {28},
	pages   = {477--539},
	doi     = {10.1146/annurev.fl.28.010196.002401}
}

@article{Dusek1994,
	author  = {J. Du\v{s}ek and P. Le Gal and P. Frauni\'e},
	title   = {A numerical and theoretical study of the first {Hopf} bifurcation in a cylinder wake},
	journal = {J. Fluid Mech.},
	year    = {1994},
	volume  = {264},
	pages   = {59--80},
	doi     = {10.1017/S0022112094000583}
}

@article{Wang2011,
	author  = {F. Wang and B. S. Lazarov and O. Sigmund},
	title   = {On projection methods, convergence and robust formulations in topology optimization},
	journal = {Struct. Multidiscip. Optim.},
	year    = {2011},
	volume  = {43},
	number  = {6},
	pages   = {767--784},
	doi     = {10.1007/s00158-010-0602-y}
}

@article{Svanberg1987,
	author  = {K. Svanberg},
	title   = {The method of moving asymptotes---a new method for structural optimization},
	journal = {Int. J. Numer. Methods Eng.},
	year    = {1987},
	volume  = {24},
	number  = {2},
	pages   = {359--373},
	doi     = {10.1002/nme.1620240207}
}

@article{Armaly1983,
	author  = {B. F. Armaly and F. Durst and J. C. F. Pereira and B. Sch{\"o}nung},
	title   = {Experimental and Theoretical Investigation of Backward-Facing Step Flow},
	journal = {J. Fluid Mech.},
	year    = {1983},
	volume  = {127},
	pages   = {473--496},
	doi     = {10.1017/S0022112083002839}
}

@article{Zhang2013,
	author  = {L. Zhang and A. Poth{\'e}rat},
	title   = {Influence of the geometry on the two- and three-dimensional dynamics of the flow in a 180$^\circ$ sharp bend},
	journal = {Phys. Fluids},
	volume  = {25},
	number  = {5},
	pages   = {053605},
	year    = {2013},
	doi     = {10.1063/1.4807070}
}

@article{Katsumata2025,
	author  = {Katsumata, R. and Nishiguchi, K. and Hoshiba, H. and Kato, J.},
	title   = {Topology Optimization for Large-Scale Unsteady Flow With the Building-Cube Method},
	journal = {Int. J. Numer. Methods Eng.},
	year    = {2025},
	volume  = {126},
	number  = {5},
	pages   = {e70016},
	doi     = {10.1002/nme.70016}
}

@article{Yaji2014,
	author  = {Yaji, K. and Yamada, T. and Yoshino, M. and Matsumoto, T. and Izui, K. and Nishiwaki, S.},
	title   = {Topology Optimization Using the Lattice {Boltzmann} Method Incorporating Level Set Boundary Expressions},
	journal = {J. Comput. Phys.},
	year    = {2014},
	volume  = {274},
	pages   = {158--181},
	doi     = {10.1016/j.jcp.2014.06.004}
}

@article{Kobayashi2021,
	author  = {Kobayashi, H. and Yaji, K. and Yamasaki, S. and Fujita, K.},
	title   = {Topology Design of Two-Fluid Heat Exchange},
	journal = {Struct. Multidiscip. Optim.},
	year    = {2021},
	volume  = {63},
	number  = {2},
	pages   = {821--834},
	doi     = {10.1007/s00158-020-02736-8}
}

@article{Hoghoj2020,
	author  = {H{\o}gh{\o}j, L. C. and N{\o}rhave, D. R. and Alexandersen, J. and Sigmund, O. and Andreasen, C. S.},
	title   = {Topology Optimization of Two Fluid Heat Exchangers},
	journal = {Int. J. Heat Mass Transfer},
	year    = {2020},
	volume  = {163},
	pages   = {120543},
	doi     = {10.1016/j.ijheatmasstransfer.2020.120543}
}

@article{Koga2013,
	author  = {Koga, A. A. and Lopes, E. C. C. and Villa Nova, H. F. and de Lima, C. R. and Silva, E. C. N.},
	title   = {Development of Heat Sink Device by Using Topology Optimization},
	journal = {Int. J. Heat Mass Transfer},
	year    = {2013},
	volume  = {64},
	pages   = {759--772},
	doi     = {10.1016/j.ijheatmasstransfer.2013.05.007}
}

@article{Wang2018,
	author  = {Wang, H. and Marshall, S. D. and Arayanarakool, R. and Balasubramaniam, L. and Jin, X. and Lee, P. S. and Chen, P. C. Y.},
	title   = {Optimization of a Fluid Distribution Manifold: Mechanical Design, Numerical Studies, Fabrication, and Experimental Verification},
	journal = {Int. J. Fluid Mach. Syst.},
	year    = {2018},
	volume  = {11},
	number  = {3},
	pages   = {244--256},
	doi     = {10.5293/IJFMS.2018.11.3.244}
}

@article{Gaymann2019,
	author  = {Gaymann, A. and Montomoli, F.},
	title   = {Fluid Topology Optimization: Bio-Inspired Valves for Aircraft Engines},
	journal = {Int. J. Heat Fluid Flow},
	year    = {2019},
	volume  = {79},
	pages   = {108455},
	doi     = {10.1016/j.ijheatfluidflow.2019.108455}
}

@article{Andreasen2009,
	author  = {Andreasen, C. S. and Gersborg, A. R. and Sigmund, O.},
	title   = {Topology Optimization of Microfluidic Mixers},
	journal = {Int. J. Numer. Methods Fluids},
	year    = {2009},
	volume  = {61},
	number  = {5},
	pages   = {498--513},
	doi     = {10.1002/fld.1964}
}

@article{Li2022,
	author  = {Li, J. and Zhu, S.},
	title   = {Shape Optimization of {Navier--Stokes} Flows by a Two-Grid Method},
	journal = {Comput. Methods Appl. Mech. Engrg.},
	year    = {2022},
	volume  = {400},
	pages   = {115531},
	doi     = {10.1016/j.cma.2022.115531}
}

@article{Chen2017,
	author  = {Chen, C. and Yaji, K. and Yamada, T. and Izui, K. and Nishiwaki, S.},
	title   = {Local-in-Time Adjoint-Based Topology Optimization of Unsteady Fluid Flows Using the Lattice {Boltzmann} Method},
	journal = {Mech. Eng. J.},
	year    = {2017},
	volume  = {4},
	number  = {3},
	pages   = {17-00120},
	doi     = {10.1299/mej.17-00120}
}

@article{Putranto2026,
	author  = {Putranto, A. B. and Zuhal, L. R. and Michelis, T. and Palar, P. S.},
	title   = {A Level-Set Topology Optimization of Drag and Lift Problems Incorporating Unsteady Flow Field Information},
	journal = {Struct. Multidiscip. Optim.},
	year    = {2026},
	volume  = {69},
	number  = {5},
	pages   = {135},
	doi     = {10.1007/s00158-026-04349-z}
}

@article{Wiker2007,
	author  = {Wiker, N. and Klarbring, A. and Borrvall, T.},
	title   = {Topology Optimization of Regions of {Darcy} and {Stokes} Flow},
	journal = {Int. J. Numer. Methods Eng.},
	year    = {2007},
	volume  = {69},
	number  = {7},
	pages   = {1374--1404},
	doi     = {10.1002/nme.1811}
}

@book{Strogatz2018,
	author    = {Strogatz, S. H.},
	title     = {Nonlinear Dynamics and Chaos: With Applications to Physics, Biology, Chemistry, and Engineering},
	publisher = {CRC Press},
	edition   = {2},
	year      = {2018},
	doi       = {10.1201/9780429492563}
}

@book{Guckenheimer1983,
	author    = {Guckenheimer, J. and Holmes, P.},
	title     = {Nonlinear Oscillations, Dynamical Systems, and Bifurcations of Vector Fields},
	publisher = {Springer},
	series    = {Applied Mathematical Sciences},
	volume    = {42},
	year      = {1983},
	doi       = {10.1007/978-1-4612-1140-2}
}

@article{Mohammed2011,
	author  = {Mohammed, H. A. and Bhaskaran, G. and Shuaib, N. H. and Saidur, R.},
	title   = {Heat transfer and fluid flow characteristics in microchannels heat exchanger using nanofluids: A review},
	journal = {Renew. Sustain. Energy Rev.},
	year    = {2011},
	volume  = {15},
	number  = {3},
	pages   = {1502--1512},
	doi     = {10.1016/j.rser.2010.11.031}
}

@article{Zhang2023,
	author  = {Zhang, J. and Zou, Z. and Fu, C.},
	title   = {A Review of the Complex Flow and Heat Transfer Characteristics in Microchannels},
	journal = {Micromachines},
	year    = {2023},
	volume  = {14},
	number  = {7},
	pages   = {1451},
	doi     = {10.3390/mi14071451}
}

@article{Jimenez2018,
	author  = {Jim{\'e}nez, J.},
	title   = {Coherent structures in wall-bounded turbulence},
	journal = {J. Fluid Mech.},
	year    = {2018},
	volume  = {842},
	pages   = {P1},
	doi     = {10.1017/jfm.2018.144}
}

@article{Chen2009Manifold,
	author  = {Chen, Andrew W. and Sparrow, Ephraim M.},
	title   = {Systematic approaches for design of distribution manifolds having the same per-port outflow},
	journal = {J. Fluids Eng.},
	year    = {2009},
	volume  = {131},
	number  = {6},
	pages   = {061101},
	doi     = {10.1115/1.3111256}
}

@article{Qian1992,
	author  = {Qian, Y. H. and d'Humi{\`e}res, D. and Lallemand, P.},
	title   = {Lattice {BGK} Models for {Navier--Stokes} Equation},
	journal = {Europhys. Lett.},
	volume  = {17},
	number  = {6},
	pages   = {479--484},
	year    = {1992},
	doi     = {10.1209/0295-5075/17/6/001}
}

@article{Ito2026MeanFlow,
	author  = {Ito, Shota and Grafen, Johannes L. and Bukreev, Fedor
	and Kummerl{\"a}nder, Adrian and Krause, Mathias J.},
	title   = {Mean-Flow Adjoint Sensitivity Analysis of Unsteady Flow
	Around Porous Cylinders Using a Homogenized Lattice
	Boltzmann Method},
	journal = {arXiv preprint arXiv:2606.31707},
	year    = {2026},
	eprint  = {2606.31707},
	archivePrefix = {arXiv},
	primaryClass  = {physics.flu-dyn},
	url     = {https://arxiv.org/abs/2606.31707}
}

@article{Liu2026,
	author  = {Liu, Z. and Li, J. and Zhu, S.},
	title   = {Topology optimization for microfluidic mixers by a phase field method},
	journal = {International Journal for Numerical Methods in Fluids},
	year    = {2026},
	doi     = {10.1002/fld.70083}
}

@article{Li2026a,
	author  = {J. Li and H. Yang and S. Zhu},
	title   = {Topology Optimization in Thermal-Fluid Coupling by a Phase-Field Method},
	journal = {SIAM J. Control Optim.},
	year    = {2026},
	volume  = {64},
	number  = {4},
	pages   = {2167--2203},
	issn    = {0363-0129},
	doi     = {10.1137/24M1639919}
}

@article{Li2026b,
	author  = {J. Li and W. Li and S. Zhu},
	title   = {Shape Optimization Constrained by Time-Dependent {Stokes} Flow: Modeling, Analysis and Numerical Simulation},
	journal = {Comput. Optim. Appl.},
	year    = {2026},
	issn    = {0926-6003},
	doi     = {10.1007/s10589-026-00807-y}
}
	
\end{document}